\documentclass[english,11pt,aps,letter,superscriptaddress,nofootinbib,pre]{revtex4}
\usepackage[T1]{fontenc}
\usepackage[latin9]{inputenc}
\usepackage{xcolor}
\usepackage{pdfcolmk}
\usepackage{float}
\usepackage{amsmath}
\usepackage{amssymb}
\usepackage{graphicx}
\usepackage{esint}
\PassOptionsToPackage{normalem}{ulem}
\usepackage{ulem}

\makeatletter

\providecommand{\tabularnewline}{\\}
\floatstyle{ruled}
\newfloat{algorithm}{tbp}{loa}
\providecommand{\algorithmname}{Algorithm}
\floatname{algorithm}{\protect\algorithmname}
\providecolor{lyxadded}{rgb}{1,0,0}
\providecolor{lyxdeleted}{rgb}{0,0,1}

\@ifundefined{textcolor}{}
{%
 \definecolor{BLACK}{gray}{0}
 \definecolor{WHITE}{gray}{1}
 \definecolor{RED}{rgb}{1,0,0}
 \definecolor{GREEN}{rgb}{0,1,0}
 \definecolor{BLUE}{rgb}{0,0,1}
 \definecolor{CYAN}{cmyk}{1,0,0,0}
 \definecolor{MAGENTA}{cmyk}{0,1,0,0}
 \definecolor{YELLOW}{cmyk}{0,0,1,0}
}

\usepackage{ae,aecompl}

\newcommand{\sfrac}[2]{\mathchoice
  {\kern0em\raise.5ex\hbox{\the\scriptfont0 #1}\kern-.15em/
   \kern-.15em\lower.25ex\hbox{\the\scriptfont0 #2}}
  {\kern0em\raise.5ex\hbox{\the\scriptfont0 #1}\kern-.15em/
   \kern-.15em\lower.25ex\hbox{\the\scriptfont0 #2}}
  {\kern0em\raise.5ex\hbox{\the\scriptscriptfont0 #1}\kern-.2em/
   \kern-.15em\lower.25ex\hbox{\the\scriptscriptfont0 #2}}
  {#1\!/#2}}

\makeatletter
\DeclareMathSizes{\@xipt}{10}{6}{5}
\makeatother

\makeatother

\usepackage{babel}
\begin{document}
\global\long\def\V#1{\boldsymbol{#1}}
\global\long\def\M#1{\boldsymbol{#1}}
\global\long\def\Set#1{\mathbb{#1}}

\global\long\def\D#1{\Delta#1}
\global\long\def\d#1{\delta#1}

\global\long\def\norm#1{\left\Vert #1\right\Vert }
\global\long\def\abs#1{\left|#1\right|}

\global\long\def\grad{\M{\nabla}}
\global\long\def\avv#1{\langle#1\rangle}
\global\long\def\av#1{\left\langle #1\right\rangle }

\global\long\def\myhalf{\sfrac{1}{2}}
\global\long\def\mythreehalves{\sfrac{3}{2}}

\global\long\def\vb{\V v}

\title{Low Mach Number Fluctuating Hydrodynamics of Binary Liquid Mixtures}

\author{Andy Nonaka}

\affiliation{Center for Computational Science and Engineering, Lawrence Berkeley
National Laboratory, Berkeley, CA, 94720}

\author{Yifei Sun}

\affiliation{Courant Institute of Mathematical Sciences, New York University,
New York, NY 10012}

\author{John B. Bell}

\affiliation{Center for Computational Science and Engineering, Lawrence Berkeley
National Laboratory, Berkeley, CA, 94720}

\author{Aleksandar Donev}

\email{donev@courant.nyu.edu}

\affiliation{Courant Institute of Mathematical Sciences, New York University,
New York, NY 10012}
\begin{abstract}
Continuing on our previous work {[}A. Donev, A. Nonaka, Y. Sun, T.
G. Fai, A. L. Garcia and J. B. Bell, Comm. App. Math. and Comp. Sci.,
9-1:47-105, 2014{]}, we develop semi-implicit numerical methods for
solving low Mach number fluctuating hydrodynamic equations appropriate
for modeling diffusive mixing in isothermal mixtures of fluids with
different densities and transport coefficients. We treat viscous dissipation
implicitly using a recently-developed variable-coefficient Stokes
solver {[}M. Cai, A. J. Nonaka, J. B. Bell, B. E. Griffith and A.
Donev, Commun. Comput. Phys., 16(5):1263-1297, 2014{]}. This allows
us to increase the time step size significantly for low Reynolds number
flows with large Schmidt numbers compared to our earlier explicit
temporal integrator. Also, unlike most existing deterministic methods
for low Mach number equations, our methods do not use a fractional
time-step approach in the spirit of projection methods, thus avoiding
splitting errors and giving full second-order deterministic accuracy
even in the presence of boundaries for a broad range of Reynolds numbers
including steady Stokes flow. We incorporate the Stokes solver into
two time-advancement schemes, where the first is suitable for inertial
flows and the second is suitable for the overdamped limit (viscous-dominated
flows), in which inertia vanishes and the fluid motion can be described
by a steady Stokes equation. We also describe how to incorporate advanced
higher-order Godunov advection schemes in the numerical method, allowing
for the treatment of (very) large Peclet number flows with a vanishing
mass diffusion coefficient. We incorporate thermal fluctuations in
the description in both the inertial and overdamped regimes. We validate
our algorithm with a series of stochastic and deterministic tests.
Finally, we apply our algorithms to model the development of giant
concentration fluctuations during the diffusive mixing of water and
glycerol, and compare numerical results with experimental measurements.
We find good agreement between the two, and observe propagative (non-diffusive)
modes at small wavenumbers (large spatial scales), not reported in
published experimental measurements of concentration fluctuations
in fluid mixtures. Our work forms the foundation for developing low
Mach number fluctuating hydrodynamics methods for miscible multi-species
mixtures of chemically reacting fluids.
\end{abstract}
\maketitle
\setlength{\abovedisplayskip}{0.35ex}\setlength{\belowdisplayskip}{0.35ex}

\section{Introduction}

Flows of realistic mixtures of miscible fluids exhibit several features
that make them more difficult to simulate numerically than flows of
simple fluids. Firstly, the physical properties of the mixture depend
on the concentration of the different species composing the mixture.
This includes both the density of the mixture at constant pressure,
and transport coefficients such as viscosity and mass diffusion coefficients.
Common simplifying assumptions such as the Boussinesq approximation,
which assumes a constant density and thus incompressible flow, or
assuming constant transport coefficients, are uncontrolled and not
appropriate for certain mixtures of very dissimilar fluids. Secondly,
for liquid mixtures there is a large separation of time scales between
the various dissipative processes, notably, mass diffusion is much
slower than momentum diffusion. The large Schmidt numbers $\mbox{Sc}\sim10^{3}-10^{4}$
typical of liquid mixtures lead to extreme stiffness and make direct
temporal integration of the hydrodynamic equations infeasible. Lastly,
flows of mixtures exhibit all of the numerical difficulties found
in single component flows, for example, well-known difficulties caused
by advection in the absence of sufficiently strong dissipation (diffusion
of momentum or mass), and challenges in incorporating thermal fluctuations
in the description. Here we develop a low Mach number approach to
isothermal binary fluid mixtures that resolves many of the above difficulties,
and paves the way for incorporating additional physics such as the
presence of more than two species \cite{MultispeciesCompressible},
chemical reactions \cite{MultispeciesChemistry,DiffusiveInstability_Chemistry_PRL},
multiple phases and surface tension \cite{StagerredFluct_Inhomogeneous,CHN_Compressible},
and others.

Stochastic fluctuations are intrinsic to fluid dynamics because fluids
are composed of molecules whose positions and velocities are random.
Thermal fluctuations affect flows from microscopic to macroscopic
scales \cite{DiffusionRenormalization_PRL,FractalDiffusion_Microgravity}
and need to be consistently included in all levels of description.
Fluctuating hydrodynamics (FHD) incorporates thermal fluctuations
into the usual Navier-Stokes-Fourier laws in the form of stochastic
contributions to the dissipative momentum, heat, and mass fluxes \cite{FluctHydroNonEq_Book}.
FHD has proven to be a very useful tool in understanding complex fluid
flows far from equilibrium \cite{LB_SoftMatter_Review,StagerredFluct_Inhomogeneous,LB_OrderParameters,SELM};
however, theoretical calculations are often only feasible after making
many uncontrolled approximations \cite{FluctHydroNonEq_Book}, and
numerical schemes used for fluctuating hydrodynamics are usually far
behind state-of-the-art deterministic computational fluid dynamics
(CFD) solvers. 

In this work, we consider binary mixtures and restrict our attention
to isothermal flows. We consider a specific equation of state (EOS)
suitable for mixtures of incompressible liquids or ideal gases, but
otherwise account for advective and diffusive mass and momentum transport
in full generality. Recently, some of us developed finite-volume methods
for the incompressible equations \cite{LLNS_Staggered}. We have also
developed low Mach number isothermal fluctuating equations \cite{LowMachExplicit},
which eliminate the stiffness arising from the separation of scales
between acoustic and vortical modes \cite{IncompressibleLimit_Majda,ZeroMach_Buoyancy,ZeroMachCombustion}.
The low Mach number equations account for the fact that for mixtures
of fluids with different densities, diffusive and stochastic mass
fluxes create local expansion and contraction of the fluid. In these
equations the incompressibility constraint should be replaced by a
``quasi-incompressibility'' constraint \cite{ZeroMachCombustion,Cahn-Hilliard_QuasiIncomp},
which introduces some difficulties in constructing conservative finite-volume
techniques \cite{LowMachAdaptive,ZeroMach_Klein,LaminarFlowChemistry,LowMach_FiniteDifference,LowMachAcoustics,LowMachExplicit}.
In Section \ref{sec:Equations} we review the low Mach number equations
of fluctuating hydrodynamics for a binary mixture of miscible fluids,
as first proposed in Ref. \cite{LowMachExplicit}.

The numerical method developed in Ref. \cite{LowMachExplicit} uses
an explicit temporal integrator. This requires using a small time
step and is infeasible for liquid mixtures due to the stiffness caused
by the separation of time scales between fast momentum diffusion and
slow mass diffusion. In recent work \cite{MultiscaleIntegrators},
some of us developed temporal integrators for the equations of fluctuating
hydrodynamics that have several important advantages. Notably, these
integrators are semi-implicit, allowing one to treat fast momentum
diffusion (viscous dissipation) implicitly, and other transport processes
explicitly. Furthermore, these temporal integrators are constructed
to be second-order accurate for the equations of linearized fluctuating
hydrodynamics (LFHD), which are suitable for describing thermal fluctuations
around stable macroscopic flows over a broad range of length and time
scales \cite{FluctHydroNonEq_Book}. Importantly, the linearization
of the fluctuating equations is carried out \emph{automatically} by
the code, making the numerical methods very similar to standard deterministic
CFD schemes. Finally, specific integrators are proposed in Ref. \cite{MultiscaleIntegrators}
to handle the extreme separation of scales between the fast velocity
and the slow concentration by taking an \emph{overdamped} limit of
the inertial equations. In this work, we extend the semi-implicit
temporal integrators proposed in Ref. \cite{MultiscaleIntegrators}
for incompressible flows to account for the quasi-incompressible nature
of low Mach number flows. We apply these temporal integrators to the
staggered-grid conservative finite-volume spatial discretization developed
in Ref. \cite{LowMachExplicit}, and additionally generalize the treatment
of advection to allow for the use of monotonicity-preserving higher-order
Godunov schemes \cite{bellColellaGlaz:1989,BDS,SemiLagrangianAdvection_2D,SemiLagrangianAdvection_3D}.

Our work relies heavily on several prior works, which we will only
briefly summarize in the present paper. The spatial discretization
we describe in more detail in Section \ref{sub:Advection} is identical
to that proposed by Donev \emph{et al} \cite{LowMachExplicit}, which
itself relies heavily on the treatment of thermal fluctuations developed
in Refs. \cite{LLNS_Staggered,LowMachExplicit}. A key development
that makes the algorithm presented here feasible for large-scale problems
is recent work by some of us \cite{StokesKrylov} on efficient multigrid-based
iterative methods for solving unsteady and steady variable-coefficient
Stokes problems on staggered grids. Our high-order Godunov method
for mass advection is based on the work of Bell \emph{et al.} \cite{BDS,SemiLagrangianAdvection_2D,SemiLagrangianAdvection_3D}.

The temporal integrators developed in Section \ref{sub:TemporalDiscretization}
are a novel approach to low Mach number hydrodynamics even in the
deterministic context. In high-resolution finite-volume methods, the
dominant paradigm has been to use a splitting (fractional-step) or
projection method \cite{Chorin68} to separate the pressure and velocity
updates \cite{almgren-iamr,LaminarFlowChemistry,LowMach_Nuclear,LowMachSupernovae,ZeroMach_Klein_2}.
We followed such a projection approach to construct an explicit temporal
integrator for the low Mach number equations \cite{LowMachExplicit}.
When viscosity is treated implicitly, however, the splitting introduces
a commutator error that leads to the appearance of spurious or ``parasitic''
modes in the presence of physical boundaries \cite{GaugeIncompressible_E,ProjectionMethods_Minion,DFDB}.
There are several techniques to reduce (but not eliminate) these artificial
boundary layers \cite{ProjectionMethods_Minion}, and for sufficiently
large Reynolds number flows the time step size dictated by advective
stability constraints makes the splitting error relatively small in
practice. At small Reynolds numbers, however, the splitting error
becomes larger as viscous effects become more dominant, and projection
methods do not apply in the steady Stokes regime for problems with
physical boundary conditions. Methods that do not split the velocity
and pressure updates but rather solve a combined Stokes system for
velocity and pressure have been used in the finite-element literature
for some time, and have more recently been used in the finite-volume
context for incompressible flow \cite{NonProjection_Griffith}. Here
we demonstrate how the same approach can be effectively applied to
the low Mach number equations for a binary fluid mixture \cite{LowMachExplicit},
to construct a method that is second-order accurate up to boundaries,
for a broad range of Reynolds numbers including steady Stokes flow. 

We test our ability to accurately capture the static structure factor
for equilibrium fluctuation calculations. Then, we test our methods
deterministically on two variable density and variable viscosity low
Mach number flows. First, we confirm second-order deterministic accuracy
in both space and time for a lid-driven cavity problem in the presence
of a bubble of a denser miscible fluid. Next, we simulate the development
of a Kevin-Helmholtz instability as a lighter less viscous fluid streams
over a denser more viscous fluid. These tests confirm the robustness
and accuracy of the methods in the presence of large contrasts, sharp
gradients, and boundaries. Next we focus on the use of fluctuating
low Mach number equations to study giant concentration fluctuations.
In Section \ref{sec:GiantFluct} we apply our methods to study the
development of giant fluctuations \cite{GiantFluctuations_Theory,GiantFluctuations_Cannell,FractalDiffusion_Microgravity,GiantFluctuations_Summary}
during free diffusive mixing of water and glycerol. We compare simulation
results to experimental measurements of the time-correlation function
of concentration fluctuations during the diffusive mixing of water
and glycerol \cite{GiantFluctuations_Cannell}. The relaxation times
show signatures of the rich deterministic dynamics, and a transition
from purely diffusive relaxation of concentration fluctuations at
large wavenumbers, to more complex buoyancy-driven dynamics at smaller
wavenumbers. We find reasonably-good agreement given the large experimental
uncertainties, and observe the appearance of propagative modes at
small wavenumbers, which we suggest could be observed in experiments
as well.

\section{\label{sec:Equations}Low Mach Number Equations}

At mesoscopic scales, in typical liquids, sound waves are very low
amplitude and much faster than momentum diffusion; hence, they can
usually be eliminated from the fluid dynamics description. Formally,
this corresponds to taking the zero Mach number singular limit $\mbox{Ma}\rightarrow0$
of the well-known compressible fluctuating hydrodynamics equations
system \cite{Landau:Fluid,FluctHydroNonEq_Book}. In the compressible
equations, the coupling between momentum and mass transport is captured
by the equation of state (EOS) for the pressure $P(\rho,c;T_{0})$
as a local function of the density $\rho(\V r,t)$ and mass concentration
$c(\V r,t)$ at a specified temperature $T_{0}\left(\V r\right)$,
assumed to be time-independent in our isothermal model. 

The low Mach number equations can be obtained by making the ansatz
that the thermodynamic behavior of the system is captured by a reference
pressure $P_{0}\left(\V r,t\right)$, with the additional pressure
contribution $\pi\left(\V r,t\right)=O\left(\mbox{Ma}^{2}\right)$
capturing the mechanical behavior while not affecting the thermodynamics.
We will restrict consideration to cases where stratification due to
gravity causes negligible changes in the thermodynamic state across
the domain. In this case, the reference pressure is spatially constant
and constrains the system so that the evolution of $\rho$ and $c$
remains consistent with the thermodynamic EOS
\begin{equation}
P\left(\rho\left(\V r,t\right),c\left(\V r,t\right);T_{0}\left(\V r\right)\right)=P_{0}\left(t\right).\label{eq:general_EOS}
\end{equation}
Physically this means that any change in concentration must be accompanied
by a corresponding change in density, as would be observed in a system
at thermodynamic equilibrium held at the fixed reference pressure
and temperature. The EOS defines density $\rho\left(c\left(\V r,t\right);T_{0}\left(\V r\right),P_{0}\left(t\right)\right)$
as an implicit function of concentration in a binary liquid mixture.
The EOS constraint (\ref{eq:general_EOS}) can be re-written as a
constraint on the divergence of the fluid velocity $\V v(\V r,t)$,
\begin{equation}
\rho\grad\cdot\V v=-\beta\,\grad\cdot\V F,\label{eq:div_v}
\end{equation}
where $\V F$ is the total diffusive mass flux defined in (\ref{eq:F_def}),
and the solutal expansion coefficient 
\[
\beta\left(c\right)=\frac{1}{\rho}\left(\frac{\partial\rho}{\partial c}\right)_{P_{0},T_{0}}
\]
is determined by the specific form of the EOS. 

In this work we consider a specific \emph{linear} EOS, 
\begin{equation}
\frac{\rho_{1}}{\bar{\rho}_{1}}+\frac{\rho_{2}}{\bar{\rho}_{2}}=\frac{c\rho}{\bar{\rho}_{1}}+\frac{(1-c)\rho}{\bar{\rho}_{2}}=1,\label{eq:EOS_quasi_incomp}
\end{equation}
where $\bar{\rho}_{1}$ and $\bar{\rho}_{2}$ are the densities of
the pure component fluids ($c=1$ and $c=0$, respectively), giving
\begin{equation}
\beta=\rho\left(\frac{1}{\bar{\rho}_{2}}-\frac{1}{\bar{\rho}_{1}}\right)=\frac{\bar{\rho}_{1}-\bar{\rho}_{2}}{c\bar{\rho}_{2}+(1-c)\bar{\rho}_{1}}.\label{eq:beta_simple}
\end{equation}
It is important that for this specific form of the EOS $\beta/\rho$
is a material constant independent of the concentration; this allows
us to write the EOS constraint (\ref{eq:div_v_constraint}) in conservative
form $\grad\cdot\V v=-\grad\cdot\left(\beta\rho^{-1}\V F\right)$
and take the reference pressure $P_{0}$ to be independent of time.
The specific form of the density dependence (\ref{eq:beta_simple})
on concentration arises if one assumes that two incompressible fluids
do not change volume upon mixing, which is a reasonable assumption
for liquids that are not too dissimilar at the molecular level. Surprisingly
the EOS (\ref{eq:EOS_quasi_incomp}) is also valid for a mixture of
ideal gases. If the specific EOS (\ref{eq:EOS_quasi_incomp}) is not
a very good approximation over the entire range of concentration $0\leq c\leq1$,
(\ref{eq:EOS_quasi_incomp}) may be a very good approximation over
the range of concentrations of interest if $\bar{\rho}_{1}$ and $\bar{\rho}_{2}$
are adjusted accordingly. Our choice of the specific form of the EOS
will aid significantly in the construction of simple conservative
spatial discretizations that strictly maintain the EOS without requiring
complicated nonlinear iterative corrections.

In fluctuating hydrodynamics, stochastic contributions to the momentum
and mass fluxes that are formally modeled as \cite{LLNS_Staggered}
\begin{align}
\M{\Sigma}=\sqrt{\eta k_{B}T}\left(\M{\mathcal{W}}+\M{\mathcal{W}}^{T}\right)\mbox{ and } & \M{\Psi}=\sqrt{2\chi\rho\mu_{c}^{-1}k_{B}T}\;\widetilde{\M{\mathcal{W}}},\label{stoch_flux_covariance-1}
\end{align}
where $k_{B}$ is Boltzmann's constant, $\eta$ is the shear viscosity,
$\chi$ is the diffusion coefficient, $\mu\left(c;T_{0},P_{0}\right)$
is the chemical potential of the mixture with $\mu_{c}=\left(\partial\mu/\partial c\right)_{P_{0},T_{0}}$,
and $\M{\mathcal{W}}(\V r,t)$ and $\widetilde{\M{\mathcal{W}}}(\V r,t)$
are standard zero mean, unit variance random Gaussian tensor and vector
fields, respectively, with uncorrelated components,
\[
\avv{\mathcal{W}_{ij}(\V r,t)\mathcal{W}_{kl}(\V r^{\prime},t')}=\delta_{ik}\delta_{jl}\;\delta(t-t^{\prime})\delta(\V r-\V r^{\prime}),
\]
and similarly for $\widetilde{\M{\mathcal{W}}}$.

A standard asymptotic low Mach analysis \cite{IncompressibleLimit_Majda},
formally treating the stochastic forcing as smooth, leads to the \emph{isothermal
low Mach number} equations for a binary mixture of fluids in conservation
form \cite{LowMachExplicit},
\begin{align}
\partial_{t}\left(\rho\V v\right)+\nabla\pi= & -\grad\cdot\left(\rho\V v\V v^{T}\right)+\grad\cdot\left(\eta\bar{\grad}\V v+\M{\Sigma}\right)+\rho\V g\label{eq:momentum_eq}\\
\partial_{t}\left(\rho_{1}\right)= & -\grad\cdot\left(\rho_{1}\V v\right)+\grad\cdot\V F\label{eq:rho1_eq}\\
\partial_{t}\left(\rho_{2}\right)= & -\grad\cdot\left(\rho_{2}\V v\right)-\grad\cdot\V F\label{eq:rho2_eq}\\
\grad\cdot\V v= & -\grad\cdot\left(\beta\rho^{-1}\V F\right),\label{eq:div_v_constraint}
\end{align}
where the deterministic and stochastic diffusive mass fluxes are denoted
by
\begin{equation}
\V F=\rho\chi\grad c+\M{\Psi}.\label{eq:F_def}
\end{equation}
Here $\bar{\grad}=\grad+\grad^{T}$ is a symmetric gradient, $\rho_{1}=\rho c$
is the density of the first component, $\rho_{2}=(1-c)\rho$ is the
density of the second component, and $\V g$ is the gravitational
acceleration. The gradient of the non-thermodynamic component of the
pressure $\pi$ (Lagrange multiplier) appears in the momentum equation
as a driving force that ensures the EOS constraint (\ref{eq:div_v_constraint})
is obeyed. We note that the bulk viscosity term gives a gradient term
that can be absorbed in $\pi$ and therefore does not explicitly need
to appear in the equations. Temperature dynamics and fluctuations
are neglected in these equations; however, this type of approach can
be extended to include thermal effects. The shear viscosity $\eta\left(c;T_{0},P_{0}\right)$
and the mass diffusion coefficient $\chi\left(c;T_{0},P_{0}\right)$
in general depend on the concentration. Note that the two density
equations (\ref{eq:rho1_eq},\ref{eq:rho2_eq}) can be combined to
obtain the usual continuity equation for the total density,
\begin{equation}
\partial_{t}\rho=-\grad\cdot\left(\rho\V v\right),\label{eq:rho_eq}
\end{equation}
and the primitive (non-conservation law) form of the concentration
equation,
\begin{equation}
\rho\left(\partial_{t}c+\V v\cdot\grad c\right)=\grad\cdot\V F.\label{eq:c_t_primitive}
\end{equation}
Our conservative numerical scheme is based on Eqs. (\ref{eq:momentum_eq},\ref{eq:rho1_eq},\ref{eq:div_v_constraint},\ref{eq:rho_eq}).

In Ref. \cite{LowMachExplicit} we discussed the effect of the low
Mach constraint on the thermal fluctuations, suitable boundary conditions
for the low Mach equations, and presented a gauge formulation of the
equations that formally eliminates pressure in a manner similar to
the projection operator formulation for incompressible flows. Importantly,
the gauge formulation demonstrates that although the low Mach equations
have the appearance of a constrained system, one can write them in
an unconstrained form by introducing a gauge degree of freedom for
the pressure. For the purposes of time integration, one can therefore
treat these equations as standard initial-value problems \cite{LowMachExplicit}
and use the temporal integrators developed in Ref. \cite{MultiscaleIntegrators}.

\subsection{\label{sub:LFHD}Linearized low Mach fluctuating hydrodynamics}

It is important to note that the equations of fluctuating hydrodynamics
should be interpreted as a mesoscopic coarse-grained representation
of the mass, momentum and energy transport in fluids \cite{OttingerBook}.
As such, these equations implicitly contain a mesoscopic coarse-graining
length and time scale that is larger than molecular scales \cite{DiscreteLLNS_Espanol}
and can only formally be written as stochastic partial differential
equations (SPDEs). A coarse-graining scale can explicitly be included
in the SPDEs \cite{DiffusionJSTAT,DDFT_Hydro}; such a coarse-graining
scale explicitly enters in our finite-volume spatio-temporal discretization
through the grid spacing (equivalently, the volume of a grid cell,
or more precisely, the number of molecules per grid cell). Additional
difficulties are posed by the fact that in general the noise in the
nonlinear equations is multiplicative, requiring a careful stochastic
interpretation; the Mori-Zwanzig projection formalism \cite{GrabertBook}
suggests the correct stochastic interpretation is the kinetic one
\cite{KineticStochasticIntegral_Ottinger}.

For compressible and incompressible flows, the SPDEs of \emph{linearized}
fluctuating hydrodynamics (LFHD) \cite{FluctHydroNonEq_Book} can
be given a precise continuum meaning \cite{DaPratoBook,LLNS_S_k,DFDB,MultiscaleIntegrators}.
In these linearized equations one splits each variable into a deterministic
component and small fluctuations around the deterministic solution,
e.g., $c\left(\V r,t\right)=\bar{c}\left(\V r,t\right)+\d c\left(\V r,t\right)$,
where $\bar{c}$ is a solution of the deterministic equations (\ref{eq:momentum_eq},\ref{eq:rho1_eq},\ref{eq:div_v_constraint},\ref{eq:rho_eq})
with $\M{\Psi}=0$ and $\M{\Sigma}=0$. Here $\d c$ is the solution
of a \emph{linear} additive-noise equation obtained by linearizing
(\ref{eq:c_t_primitive}) to first order in the fluctuations and evaluating
the noise amplitude at the deterministic solution; more precisely,
LFHD is an expression of the central limit theorem in the limit of
weak noise. In this work, in the stochastic setting we restrict our
attention to LFHD equations. As discussed in Ref. \cite{MultiscaleIntegrators},
we do not need to write down the (rather tedious) complete form of
the linearized low Mach number equations (for an illustration, see
Section \ref{sub:Overdamped}) since the numerical method will perform
this linearization automatically. Namely, the complete nonlinear equations
are essentially equivalent to the LFHD equations when the noise is
sufficiently weak, i.e., when the hydrodynamic cells contain many
molecules.

The low Mach number equations pose additional difficulties because
they represent a coarse-graining of the dynamics not just in space
but also in time. As such, even the linearized equations cannot directly
be interpreted as describing a standard diffusion process. This is
because the stochastic mass flux $\M{\Psi}$ in the EOS constraint
(\ref{eq:div_v_constraint}) makes the velocity formally white-in-time
\cite{LowMachExplicit}. We note, however, that the analysis in Ref.
\cite{DiffusionJSTAT} shows that there is a close connection between
mass diffusion and advection by the thermally-fluctuating velocity
field, and thus between $\M{\Psi}$ and velocity fluctuations. This
suggests that a precise interpretation of the low Mach constraint
in the presence of stochastic mass fluxes requires a very delicate
mathematical analysis. In this work we rely on the implicit coarse-graining
in time provided by the finite time step size in the temporal integration
schemes to regularize the low Mach equations \cite{LowMachExplicit}.
Furthermore, for the applications we study here, we can neglect stochastic
mass fluxes and assume $\M{\Psi}\approx0$, in which case the difficulties
related to a white-in-time velocity disappear.

\subsection{\label{sub:Overdamped}Overdamped limit}

At small scales, flows in liquids are viscous-dominated and the inertial
momentum flux $\rho\V v\V v^{T}$ can often be neglected in a zero
Reynolds number approximation. In addition, in liquids, there is a
large separation of time scales between the fast momentum diffusion
and slow mass diffusion, i.e., the Schmidt number $\mbox{Sc}=\eta/\left(\rho\chi\right)$
is large. This makes the relaxation times of velocity modes at sufficiently
large wavenumbers much smaller than those of the concentration modes.
Formally treating the stochastic force terms as smooth for the moment,
the separation of time scales implies that we can replace the \emph{inertial}
momentum equation (\ref{eq:momentum_eq}) with the \emph{overdamped}
steady-Stokes equation

\begin{eqnarray}
-\grad\cdot\left(\eta\bar{\grad}\V v\right)+\nabla\pi & = & \grad\cdot\M{\Sigma}+\rho\V g\label{eq:steady_Stokes}\\
\grad\cdot\V v & = & -\left(\bar{\rho}_{2}^{-1}-\bar{\rho}_{1}^{-1}\right)\grad\cdot\left(\chi\rho\grad c+\M{\Psi}\right).\nonumber 
\end{eqnarray}
The above equations can be used to eliminate velocity as a variable,
leaving only the concentration equation (\ref{eq:c_t_primitive}).
Note that the density equation (\ref{eq:rho_eq}) simply defines density
as a function of concentration and thus is not considered an independent
equation. 

The solution of the Stokes system
\begin{eqnarray}
-\grad\cdot\left(\eta\bar{\grad}\V v\right)+\nabla\pi & = & \V f\label{eq:steady_Stokes_general}\\
\grad\cdot\V v & = & -h,\nonumber 
\end{eqnarray}
where $\V f\left(\V r,t\right)$ and $h\left(\V r,t\right)$ are applied
forcing terms, can be expressed in terms of a generalized inverse
Stokes linear operator %
\footnote{More generally, in the presence of inhomogeneous boundary conditions,
the solution operator for (\ref{eq:steady_Stokes_general}) is an
affine rather than a linear operator.%
} $\M{\mathcal{L}}^{-1}\left[\eta\left(\cdot,t\right)\right]$ that
is a \emph{functional} of the viscosity (and thus the concentration),
\[
\V v=\M{\mathcal{L}}^{-1}\left[\eta\right]\left(\V f,\, h\right).
\]

In the linearized fluctuating equations, one must linearize around
the (time dependent) solution of the deterministic nonlinear equation
\begin{equation}
\bar{\rho}\left(\partial_{t}\bar{c}+\bar{\V v}\cdot\grad\bar{c}\right)=\grad\cdot\left(\bar{\rho}\bar{\chi}\grad\bar{c}\right),\label{eq:c_bar_overdamped}
\end{equation}
 where we have used the shorthand notation $\bar{\rho}=\rho\left(\bar{c}\right)$,
$\bar{\eta}=\eta\left(\bar{c}\right)$, $\bar{\chi}=\chi\left(\bar{c}\right)$.
Here the velocity is an implicit function of concentration defined
via
\begin{eqnarray*}
-\grad\cdot\left(\bar{\eta}\bar{\grad}\V v\right)+\nabla\bar{\pi} & = & \bar{\rho}\V g\\
\grad\cdot\bar{\V v} & = & -\left(\bar{\rho}_{2}^{-1}-\bar{\rho}_{1}^{-1}\right)\grad\cdot\left(\bar{\chi}\bar{\rho}\grad\bar{c}\right),
\end{eqnarray*}
which we can write in shorthand notation as
\begin{equation}
\bar{\V v}=\M{\mathcal{L}}^{-1}\left[\bar{\eta}\right]\left(\bar{\rho}\V g,\;\left(\bar{\rho}_{2}^{-1}-\bar{\rho}_{1}^{-1}\right)\grad\cdot\left(\bar{\chi}\bar{\rho}\grad\bar{c}\right)\right).\label{eq:v_bar_overdamped}
\end{equation}
Here we develop second-order integrators for the deterministic overdamped
low Mach equation (\ref{eq:c_bar_overdamped},\ref{eq:v_bar_overdamped}).

In the stochastic setting, the solution of (\ref{eq:steady_Stokes})
is white in time because the stochastic mass and momentum fluxes are
white in time. This means that the advective term $\V v\cdot\grad c$
requires a specific stochastic interpretation, in addition to the
usual regularization (smoothing) in space required to interpet all
nonlinear terms appearing in formal fluctuating hydrodynamics SPDEs.
By performing a precise (albeit formal) adiabatic mode elimination
of the fast velocity variable under the assumption of infinite separation
of time scales, Donev \emph{et al.} arrive at a Stratonovich interpretation
of the random advection term $\V v\cdot\grad c$ (see Appendix A of
Ref. \cite{DiffusionJSTAT}). This analysis does not, however, directly
extend to the low Mach number equations since it relies in key ways
on the incompressibility of the fluid. Generalizing this sort of analysis
to the case of variable fluid density is nontrivial, likely requiring
the use of the gauge formulation of the low Mach equations, and appears
to be beyond the scope of existing techniques. Variable (i.e., concentration-dependent)
viscosity and mass diffusion coefficient can be handled using existing
techniques although there are subtle nonlinear stochastic effects
arising from the fact that the noise in the velocity equation is multiplicative
and the invariant measure (equilibrium distribution) of the fast velocity
depends on the slow concentration.

In the linearized setting, however, the difficulties associated with
the interpretation of stochastic integrals and multiplicative noise
disappear. The complete form of the linearized equations contains
many terms and is rather tedious. Since we will never need to explicitly
write this form let us illustrate the procedure by assuming $\chi$
and $\eta$ to be constant. For the concentration, we obtain the linearized
equation
\begin{equation}
\bar{\rho}\left(\partial_{t}\left(\d c\right)+\left(\d{\V v}\right)\cdot\grad\bar{c}\right)=\grad\cdot\left(\bar{\rho}\chi\grad\left(\d c\right)+\bar{\rho}^{\prime}\chi\left(\grad\bar{c}\right)\d c\right)-\bar{\rho}^{-1}\bar{\rho}^{\prime}\grad\cdot\left(\bar{\rho}\chi\grad\bar{c}\right)\d c,\label{eq:c_t_lin}
\end{equation}
where $\bar{\rho}^{\prime}=d\rho\left(\bar{c}\right)/dc=\bar{\rho}\beta\left(\bar{\rho}\right)$
relates concentration fluctuations to density fluctuations via the
EOS. Here we split $\d{\V v}=\d{\V v}_{c}+\d{\V v}_{f}$ into a component
$\d{\V v}_{c}$ that is continuous in time and a component $\d{\V v}_{f}$
that is white in time,
\begin{eqnarray*}
\d{\V v}_{c} & = & \M{\mathcal{L}}^{-1}\left[\bar{\eta}\right]\left(\bar{\rho}^{\prime}\V g\d c,\;\left(\bar{\rho}_{2}^{-1}-\bar{\rho}_{1}^{-1}\right)\grad\cdot\left(\bar{\rho}\chi\grad\left(\d c\right)\right)\right),\\
\d{\V v}_{f} & = & \M{\mathcal{L}}^{-1}\left[\bar{\eta}\right]\left(\grad\cdot\M{\Sigma},\;\grad\cdot\M{\Psi}\right).
\end{eqnarray*}
The term $\bar{\rho}\left(\d{\V v}_{f}\right)\cdot\grad\bar{c}$ in
(\ref{eq:c_t_lin}) is interpreted as an additive noise term with
a rather complicated and potentially time-dependent (via $\bar{\eta}\left(\bar{c}\left(\V r,t\right)\right)$)
spatial correlation structure. In this work we develop numerical methods
that solve the overdamped linearized equation (\ref{eq:c_t_lin})
to second-order weakly \cite{MultiscaleIntegrators}.

\section{\label{sec:SpatioTemporal}Spatio-Temporal Discretization}

Our baseline spatio-temporal discretization of the low Mach equations
is based on the method of lines approach where we first discretize
the (S)PDEs in space to obtain a system of (S)ODEs, which we then
solve using a single-step multi-stage temporal integrator. The conservative
finite-volume spatial discretization that we employ here is essentially
identical to that developed in our previous works \cite{LowMachExplicit}
and \cite{StokesKrylov}. In summary, scalar fields such as concentration
and densities are cell-centered, while velocity is face-centered.
In order to ensure conservation, the conserved momentum $\rho\V v$
and mass densities $\rho_{1}$ and $\rho$ are evolved rather than
the primitive variables $\V v$ and $c$. Diffusion of mass and momentum
is discretized using standard centered differences, leading to compact
stencils similar to the standard Laplacian. Stochastic mass fluxes
are associated with the faces of the regular grid, while for stochastic
momentum fluxes we associate the diagonal elements with the cell centers
and the off-diagonal elements with the nodes (in 2D) or edges (in
3D) of a regular grid with grid spacing $\D x$.

Here we focus our discussion on three new aspects of our spatio-temporal
discretization. After summarizing the dimensionless numbers that control
the appropriate choice of advection method and temporal integrator,
in Section \ref{sub:Advection} we describe our implementation of
two advection schemes and a discussion of the advantages of each.
In Section \ref{sub:GMRES} we describe our implicit treatment of
viscous dissipation using a GMRES solver for the coupled velocity-pressure
Stokes system. In Section \ref{sub:TemporalDiscretization} we describe
our overall temporal discretization strategies for the inertial and
overdamped regimes.

\subsection{Dimensionless Numbers}

The suitability of a particular temporal integrator or advection scheme
depends on the following dimensionless numbers:
\begin{eqnarray*}
\mbox{cell Reynolds number} &  & \mbox{Re}_{c}=\frac{U\D x}{\nu}\\
\mbox{cell Peclet number} &  & \mbox{Pe}_{c}=\frac{U\D x}{\chi}\\
\mbox{Schmidt number} &  & \mbox{Sc}=\frac{\nu}{\chi}=\frac{\mbox{Pe}}{\mbox{Re}},
\end{eqnarray*}
where $\nu=\eta/\rho$ is the kinematic viscosity. Observe that the
first two depend on the spatial resolution and the typical flow speed
$U$, while the Schmidt number is an intrinsic material property of
the mixture. Also note that the physically relevant Reynolds $\mbox{Re}$
and Peclet $\mbox{Pe}$ numbers would be defined with a length scale
much larger than $\D x$, such as the system size, and thus would
be much larger than the discretization-scale numbers above. In this
work, we are primarily interested in small-scale flows with $\mbox{Re}_{c}\lesssim1$
and large $\mbox{Sc}$ (liquid mixtures).

The choice of advection scheme for concentration (partial densities)
is dictated by $\mbox{Pe}_{c}$. If $\mbox{Pe}_{c}\gtrsim2$, centered
advection schemes will generate non-physical oscillations, and one
must use the Godunov advection scheme described below. However, it
is important to note that in this case the spectrum of fluctuations
will not be correctly preserved by the advection scheme; if fluctuations
need to be resolved it is advisable to instead reduce the grid spacing
and thus reduce the cell Peclet number to $\mbox{Pe}_{c}\lesssim2$
and use centered advection.

The choice of the temporal integrator, on the other hand, is determined
by the importance of inertia and the time scale of interest. If $\mbox{Re}_{c}$
is not sufficiently small, then there is no alternative to resolving
the inertial dynamics of the velocity. Now let us assume that $Re\ll1$,
i.e., viscosity is dominant. If the time scale of interest is the
advective timescale $L/U$, where $L$ is the system size, then one
should use the inertial equations. However, the inertial temporal
integrator described in Section \ref{sub:TemporalDiscretization}
will be rather inefficient if the time scale of interest is the diffusive
time scale $L^{2}/\chi$, as is the case in the study of diffusive
mixing presented in Section \ref{sec:GiantFluct}. This is because
the Crank-Nicolson (implicit midpoint) scheme used to treat viscosity
in our methods is only A-stable, and, therefore, if the viscous Courant
number $\nu\D t/\D x^{2}$ is too large, unphysical oscillations in
the solution will appear (note that this problem is much more serious
for fluctuating hydrodynamics due to the presence of fluctuations
at \emph{all} scales). In order to be able to use a time step size
on the diffusive time scale, one must construct a \emph{stiffly accurate}
temporal integrator. This requires using an L-stable scheme to treat
viscosity, such as the backward Euler scheme, which is however only
first-order accurate.

Constructing a second-order stiffly accurate implicit-explicit integrator
in the context of variable density low Mach flows is rather nontrivial.
Furthermore, using an L-stable scheme leads to a damping of the velocity
fluctuations at large wavenumbers and is inferior to the implicit
midpoint scheme in the context of fluctuating hydrodynamics \cite{DFDB}.
Therefore, in this work we choose to consider separately the overdamped
limit $\mbox{Re}\rightarrow0$ and $\mbox{Sc}\rightarrow\infty$ (note
that the value of $\mbox{Pe}$ is arbitrary). In this limit we analytically
eliminate the velocity as an independent variable, leaving only the
concentration equation, which evolves on the diffusive time scale.
We must emphasize, however, that the overdamped equations should be
used with caution, especially in the presence of fluctuations. Notably,
the validity of the overdamped approximation requires that the separation
of time scales between the fast velocity and slow concentration be
uniformly large over \emph{all} wavenumbers, since fluctuations are
present at \emph{all} lengthscales. In the study of giant fluctuations
we present in Section \ref{sec:GiantFluct}, buoyancy effects speed
up the dynamics of large-scale concentration fluctuations and using
the overdamped limit would produce physically incorrect results at
small wavenumbers. In microgravity, however, the overdamped limit
is valid and we have used it to study giant fluctuations over very
long time scales in a number of separate works \cite{GRADFLEXTransient,DiffusionJSTAT}.

\subsection{\label{sub:Advection}Advection}

We have implemented two advection schemes for cell-centered scalar
fields, and describe under what conditions each is more suitable.
The first is a simpler non-dissipative centered advection discretization
described in our previous work \cite{LowMachExplicit}. This scheme
preserves the skew-adjoint nature of advection and thus maintains
fluctuation-dissipation balance in the stochastic context. However,
when sharp gradients are present, centered advection schemes require
a sufficient amount of dissipation (diffusion) in order to avoid the
appearance of Gibbs-phenomenon instabilities. Higher-order Godunov
schemes have been used successfully with cell-centered finite volume
schemes for some time \cite{bellColellaGlaz:1989,BDS,SemiLagrangianAdvection_2D,SemiLagrangianAdvection_3D}.
In these semi-Lagrangian advection schemes, a construction based on
characteristics is used to estimate the average value of the advected
quantity passing through each cell face during a time step. These
averages are then used to evaluate the advective fluxes. Our second
scheme for advection is the higher-order Godunov approach of Bell,
Dawson, and Schubin (BDS) \cite{BDS}. Additional details of this
approach are provided in Section \ref{sub:BDS}.

The BDS scheme can only be used to advect cell-centered scalar fields
such as densities. This is because the scheme operatse on control
volumes, and therefore applying it to staggered field requires the
use of disjoint control volumes, thereby greatly complicating the
advection procedure for non-cell-centered data. We therefore limit
ourselves to using the skew-adjoint centered advection scheme described
in Refs. \cite{LLNS_Staggered,DFDB} to advect momentum. Although
some Godunov schemes for advecting a staggered momentum field have
been developed \cite{StaggeredGodunov,NonProjection_Griffith}, they
are not at the same level of sophistication as those for cell-centered
scalar fields. For example, in Ref. \cite{StaggeredGodunov} a piecewise
constant reconstruction is used, and in Ref. \cite{NonProjection_Griffith}
only extrapolation in space is performed but not in time. In our target
applications, there is sufficient viscous dissipation to stabilize
centered advection of momentum (note that the mass diffusion coefficient
is several orders of magnitude smaller than the kinematic viscosity
in typical liquids).

The BDS advection scheme is not skew-adjoint and thus adds some dissipation
in regions of sharp gradients that are not resolved by the underlying
grid. Thus, unlike the case of using centered advection, the spatially
discrete (but still continuous in time) fluctuating equations do not
obey a strict discrete fluctuation-dissipation principle \cite{LLNS_S_k,DFDB}.
Nevertheless, in high-resolution schemes such as BDS artificial dissipation
is added locally in regions where centered advection would have failed
completely due to insufficient spatial resolution. Furthermore, the
BDS scheme offers many advantages in the deterministic context and
allows us to simulate high Peclet number flows with little to no mass
diffusion. For well-resolved flows with sufficient dissipation there
is little difference between BDS and centered advection. Note that
both advection schemes are spatially second-order accurate for smooth
flows.

\subsubsection{\label{sub:BDS}BDS advection}

Simple advection schemes, such as the centered scheme described in
our previous work \cite{LowMachExplicit}, directly computes the divergence
of the advective flux $\V f=\phi\vb$ evaluated at a specific point
in time, where $\phi$ is a cell-centered quantity such as density,
and $\V v$ is a specified face-centered velocity. By contrast, the
BDS scheme uses the multidimensional characteristic geometry of the
advection equation 
\begin{equation}
\frac{\partial\phi}{\partial t}+\nabla\cdot(\phi\vb)=q,\label{eq:adv_diff}
\end{equation}
to estimate time-averaged fluxes through cell faces over a time interval
$\D t$, given $\phi^{n}$, as well as a face-centered velocity field
$\vb$ and a cell-centered source $q$ that are assumed \emph{constant}
over the time interval. In actual temporal discretizations $\vb\approx\V v(t^{n+\myhalf})\approx\V v(t^{n}+\D t/2)$
is a midpoint (second-order) approximation of the velocity over the
time step. Similarly, $q\approx q(t^{n+\myhalf})$ will be a centered
approximation of the divergence of the diffusive and stochastic fluxes
over the time step. In the description of our temporal integrators,
we will use the shorthand notation $\mbox{BDS}$ to denote the approximation
to the advective fluxes used in the BDS scheme for solving (\ref{eq:adv_diff}),
\[
\phi^{n+1}=\phi^{n}-\D t\,\grad\cdot\left(\text{BDS}\left(\phi^{n},\,\V v,\, q,\,\D t\right)\right)+\D t\, q.
\]

BDS is a conservative scheme based on computing time-averaged advective
fluxes through every face of the computational grid, for example,
in two dimensions,
\[
\text{BDS}_{i+\myhalf,j}=f_{i+\myhalf,j}=\phi_{i+\myhalf,j}v_{i+\myhalf,j},
\]
where $v_{i+\myhalf,j}$ is the given normal velocity at the face,
and $\phi_{i+\myhalf,j}$ represents the space-time average of $\phi$
passing through face-$(i+\frac{1}{2},\; j)$ in the time interal $\Delta t$.
The extrapolated face-centered states $\phi_{i+\myhalf,j}$ are computed
by first reconstructing a piecewise continuous profile of $\phi\left(\V r,t\right)$
in every cell that can, optionally, be limited based on monotonicity
considerations. The multidimensional characteristic geometry of the
flow in space-time is then used to estimate the time-averaged flux;
see the original papers \cite{BDS,SemiLagrangianAdvection_2D,SemiLagrangianAdvection_3D}
for a detailed description. In the original advection BDS schemes
in two dimensions \cite{BDS} and three dimensions \cite{SemiLagrangianAdvection_3D},
a piecewise-bilinear (in two dimensions) or trilinear (in three dimensions)
reconstruction of $\phi$ was used. Subsequently, the schemes were
extended to a quadratic reconstruction in two dimensions \cite{SemiLagrangianAdvection_2D}.
Note that handling boundary conditions in BDS properly requires additional
investigations, and the construction of specialized one-sided reconstruction
stencils near boundaries. In our implementation we rely on cubic extrapolation
based on interior cells and the specified boundary condition (Dirichlet
or Neumann) to fill ghost cell values behind physical boundaries,
and then apply the BDS procedure to the interior cells using the extrapolated
ghost cell values.

BDS advection, as described in \cite{BDS,SemiLagrangianAdvection_2D,SemiLagrangianAdvection_3D},
does not strictly preserve the EOS constraint, unlike centered advection.
The characteristic extrapolation of densities to space-time midpoint
values on the faces of the grid, $\left(\rho_{1}\right)_{i+\myhalf,j}$
and $\rho_{i+\myhalf,j}$, are not necessarily consistent with the
EOS, unlike centered advection where they are simple averages of values
from neighboring cells, and thus guaranteed to obey the EOS by linearity.
A simple fix that makes BDS preserve the EOS, without affecting its
formal order of accuracy, is to enforce the EOS on each face by projecting
the extrapolated values $\left(\rho_{1}\right)_{i+\myhalf,j}$ and
$\rho_{i+\myhalf,j}$ onto the EOS. In the $L_{2}$ sense, such a
projection consists of the update
\[
\left(\rho_{1}\right)_{i+\myhalf,j}\leftarrow\frac{\bar{\rho}_{1}^{2}}{\bar{\rho}_{1}^{2}+\bar{\rho}_{2}^{2}}\left(\rho_{1}\right)_{i+\myhalf,j}-\frac{\bar{\rho}_{1}\bar{\rho}_{2}}{\bar{\rho}_{1}^{2}+\bar{\rho}_{2}^{2}}\left(\rho_{2}\right)_{i+\myhalf,j},
\]
and similarly for $\rho_{2}$, or equivalently, $\rho=\rho_{1}+\rho_{2}$.
Note that this projection is done on each face only for the purposes
of computing advective fluxes and is distinct from any projection
onto the EOS performed globally.

\subsection{\label{sub:GMRES}GMRES solver}

The temporal discretization described in our previous work \cite{LowMachExplicit}
was fully explicit, whereas the discretization we employ here is implicit
in the viscous dissipation. The implicit treatment of viscosity is
traditionally handled by time-splitting approaches, in which a velocity
system is solved first, without strictly enforcing the constraint.
The solution is then projected onto the space of vector fields satisfying
the constraint \cite{Chorin68}. This type of time-splitting introduces
several artifacts, especially for viscous-dominated flows; here we
avoid time-splitting by solving a combined velocity-pressure Stokes
linear system, as discussed in detail in Ref. \cite{StokesKrylov}. 

The implicit treatment of viscosity in the temporal integrators described
in Section \ref{sub:TemporalDiscretization} requires solving discretized
unsteady Stokes equations for a velocity $\V v$ and a pressure $\pi$,

\begin{eqnarray*}
\theta\rho\V v-\grad\cdot\left(\eta\bar{\grad}\V v\right)+\grad\pi & = & \V f\\
\grad\cdot\V v & = & h,
\end{eqnarray*}
for given spatially-varying density $\rho$ and viscosity $\eta$,
right-hand sides $\V f$ and $h$, and a coefficient $\theta\geq0$.
We solve these linear systems using a GMRES Krylov solver preconditioned
by the multigrid-based preconditioners described in detail in Ref.
\cite{StokesKrylov}. This approach requires only standard velocity
(Helmholtz) and pressure (Poisson) multigrid solvers, and requires
about two-three times more multigrid iterations than solving an uncoupled
pair of velocity and pressure subproblems (as required in projection-based
splitting methods).

There are two issues that arise with the Stokes solver in the context
of temporal integration that need special care. In fluctuating hydrodynamics,
typically the average flow $\bar{\V v}$ changes slowly and is much
larger in magnitude than the fluctuations around the flow $\d{\V v}$.
In the predictor stages of our temporal integrators the convergence
criterion in the GMRES solver is based on relative tolerance. Because
the right-hand side of the linear system and the residual are dominated
by the deterministic flow, it is hard to determine when the fluctuating
component of the flow has converged to the desired relative accuracy.
In the corrector stage of our predictor-corrector schemes, we use
the predicted state as a reference, and switch to using absolute error
as the convergence criterion in GMRES, using the same residual error
tolerance as was used in the predictor stage. This ensures that the
corrector stage GMRES converges quickly if the predicted state is
already a sufficiently accurate solution of the Stokes system. Another
issue that has to be handled carefully is the imposition of inhomogeneous
boundary conditions, which leads to a linear system of the form
\[
\M A\V x^{{\rm new}}+\V b_{BC}=\V b,
\]
where $\V b_{BC}$ comes from non-homogeneous boundary conditions.
Both of these problems are solved by using a residual correction technique
to convert the Stokes linear system into one for the \emph{change}
in the velocity and pressure $\D{\V x}=\V x^{{\rm new}}-\overline{\V x}^{{\rm old}}$
relative to an initial guess or \emph{reference }state $\overline{\V x}^{{\rm old}}$,
which is typically the last known velocity and pressure, \emph{except}
that the desired inhomogeneous boundary conditions are imposed; this
ensures that boundary terms vanish and the Stokes problem for $\D{\V x}$
is in homogeneous form. Note that any Dirichlet boundary conditions
for the normal component of velocity should be consistent with $h$,
and any Dirchlet boundary conditions for the tangential component
of the velocity should be evaluated at the same point in time (e.g.,
beginning, midpoint, or endpoint of the time step) as $h$.

\subsection{\label{sub:TemporalDiscretization}Temporal Discretization}

In this section we construct temporal integrators for the spatially-discretized
low Mach number equations, in which we treat viscosity semi-implicitly.
For our target applications, the Reynolds number is sufficiently small
and the Schmidt number is sufficiently large that an explicit viscosity
treatment would lead to an overall viscous time step restriction,

\[
\frac{\eta\D t}{\D x^{2}}<\frac{1}{2d},
\]
We present temporal integrators in which we avoid fractional time
stepping and ensure \emph{strict} (to within solver and roundoff tolerances)
conservation and preservation of the EOS constraint. The key feature
of the algorithms developed here is the implicit treatment of viscous
dissipation, without, however, using splitting between the velocity
and pressure updates, as discussed at more length in the introduction.
The feasibility of this approach relies an efficient solver for Stokes
systems on a staggered grid \cite{StokesKrylov}; see also Section
\ref{sub:GMRES} for additional details.

In the temporal integrators developed here, we treat advection explicitly,
which limits the advective Courant number to
\[
\frac{v_{\text{max}}\D t}{\D x}<C\sim1.
\]
Mass diffusion is also treated explicitly since it is typically much
slower than momentum diffusion and in many examples also slower than
advection. Explicit treatment of mass diffusion leads to an additional
stability limit on the time step since the diffusive Courant number
must be sufficiently small,
\[
\frac{\chi\D t}{\D x^{2}}<\frac{1}{2d},
\]
where $d$ is the number of spatial dimensions and $\D x$ is the
grid spacing. Implicit treatment of mass diffusion is straightforward
for incompressible flows, see Algorithm 2 in Ref \cite{MultiscaleIntegrators},
but is much harder for the low Mach number equations due to the need
to maintain the EOS constraint (\ref{eq:EOS_quasi_incomp}) via the
constraint (\ref{eq:div_v_constraint}). Even with explicit mass diffusion,
provided that the Reynolds is sufficiently small and the Schmidt number
is sufficiently large, a semi-implicit viscosity treatment results
in a much larger allowable time step.

\subsubsection{Predictor-Corrector Time Stepping Schemes}

In Algorithm \ref{alg:LMInertial} we give the steps involved in advancing
the solution from time level $n$ by a time interval $\D t$ to time
level $n+1$, using a semi-implicit trapezoidal temporal integrator
\cite{MultiscaleIntegrators} for the inertial fluctuating low Mach
number equations (\ref{eq:momentum_eq},\ref{eq:rho1_eq},\ref{eq:div_v_constraint},\ref{eq:rho_eq}).
In Algorithm \ref{alg:LMOverdamped} we give an explicit midpoint
temporal integrator \cite{MultiscaleIntegrators} for the overdamped
low Mach number equations (\ref{eq:rho1_eq},\ref{eq:rho_eq},\ref{eq:steady_Stokes}).

In order to ensure strict conservation of mass and momentum, we evolve
the momentum density $\M m=\rho\V v$ and the mass densities $\rho_{1}$
and $\rho$ (an equally valid choice is to evolve $\rho_{1}$ and
$\rho_{2}$). Whenever required, the primitive variables $\V v=\M m/\rho$
and $c=\rho_{1}/\rho$ are computed from the conserved quantities.
Unlike the incompressible equations, the low Mach number equations
require the enforcement of the EOS constraint (\ref{eq:EOS_quasi_incomp})
at every update of the mass densities $\rho_{1}$ and $\rho$, notably,
both in the predictor and the corrector stages. This requires that
the right-hand side of the velocity constraint (\ref{eq:div_v_constraint})
be consistent with the corresponding diffusive fluxes used to update
$\rho_{1}$. In order to preserve the EOS and also maintain strict
conservation, Algorithms \ref{alg:LMInertial} and \ref{alg:LMOverdamped}
use a splitting approach, in which we first update the mass densities
and then we update the velocity using the updated values for the density
$\rho$ and the diffusive fluxes that will be used to update $\rho_{1}$
. Note, however, that after many time steps the small errors in enforcing
the EOS due to roundoff and solver tolerances can accumulate and lead
to a systematic drift from the EOS. This can be corrected by periodically
projecting the solution back onto the EOS using an $L_{2}$ projection,
see Section III.C in Ref. \cite{LowMachExplicit}.

In our presentation of the temporal integrators, we use superscripts
to denote where a given quantity is evaluated, for example, $\eta^{p,n+1}\equiv\eta\left(c^{p,n+1}\right)$.
Even though we use continuum notation for the divergence, gradient
and Laplacian operators, it is implicitly understood that the equations
have been discretized in space. The white-noise random tensor fields
$\M{\mathcal{W}}(\V r,t)$ and $\widetilde{\M{\mathcal{W}}}(\V r,t)$
are represented via one or two collections of i.i.d. uncorrelated
normal random variables $\V W$ and $\widetilde{\V W}$, generated
independently at each time step, as indicated by superscripts and
subscripts \cite{DFDB,MultiscaleIntegrators}. Spatial discretization
adds an additional factor of $\D V^{-\myhalf}$ due to the delta function
correlation of white-noise, where $\D V$ is the volume of a grid
cell \cite{DFDB}. For simplicity of notation we denote $\overline{\M W}=\M W+\M W^{T}$.

\begin{algorithm}
\caption{\label{alg:LMInertial}Semi-implicit trapezoidal temporal integrator
for the inertial fluctuating low Mach number equations (\ref{eq:momentum_eq},\ref{eq:rho1_eq},\ref{eq:div_v_constraint},\ref{eq:rho_eq}).}

\begin{enumerate}
\item Compute the diffusive / stochastic fluxes for the predictor. Note
that these can be obtained from step \ref{enu:InertialCorrector}
of the previous time step,
\[
\V F^{n}=\left(\rho\chi\grad c\right)^{n}+\sqrt{\frac{2\left(\chi\rho\mu_{c}^{-1}\right)^{n}k_{B}T}{\D t\,\D V}}\,\widetilde{\V W}^{n}.
\]

\item Take a predictor forward Euler step for $\rho_{1}$, and similarly
for $\rho_{2}$, or, equivalently, for $\rho$,
\begin{eqnarray*}
\rho_{1}^{\star,n+1} & = & \rho_{1}^{n}+\D t\,\grad\cdot\V F^{n}-\D t\,\grad\cdot\begin{cases}
\text{BDS}\left(\rho_{1}^{n},\,\V v^{n},\,\grad\cdot\V F^{n},\,\D t\right) & \mbox{ for BDS}\\
\rho_{1}^{n}\V v^{n} & \mbox{ for centered}
\end{cases}.
\end{eqnarray*}

\item Compute $c^{\star,n+1}=\rho_{1}^{\star,n+1}/\rho^{\star,n+1}$ and
calculate corrector diffusive fluxes and stochastic fluxes, 
\[
\V F^{\star,n+1}=\left(\rho\chi\grad c\right)^{\star,n+1}+\sqrt{\frac{2\left(\chi\rho\mu_{c}^{-1}\right)^{\star,n+1}k_{B}T}{\D t\,\D V}}\,\widetilde{\V W}^{n}.
\]

\item Take a predictor Crank-Nicolson step for the velocity, using $\V v^{n}$
as a reference state for the residual correction form of the Stokes
system,
\begin{eqnarray*}
\frac{\rho^{\star,n+1}\V v^{\star,n+1}-\rho^{n}\V v^{n}}{\D t}+\grad\pi^{\star,n+1} & = & \grad\cdot\left(-\rho\V v\V v\right)^{n}+\rho^{n}\V g\\
 & + & \frac{1}{2}\grad\cdot\left[\left(\eta^{n}\bar{\grad}\V v^{n}+\eta^{\star,n+1}\bar{\grad}\V v^{\star,n+1}\right)\right]+\grad\cdot\left(\sqrt{\frac{\eta^{n}k_{B}T}{\D t\,\D V}}\,\overline{\M W}^{n}\right),\\
\grad\cdot\V v^{\star,n+1} & = & -\grad\cdot\left(\beta\rho^{-1}\V F^{\star,n+1}\right).
\end{eqnarray*}
\label{enu:CorrAdv}Take a corrector step for $\rho_{1}$, and similarly
for $\rho_{2}$, or, equivalently, for $\rho$,
\begin{eqnarray*}
\rho_{1}^{n+1} & = & \rho_{1}^{n}+\frac{\D t}{2}\,\grad\cdot\V F^{\star,n+\myhalf}-\D t\,\grad\cdot\begin{cases}
\text{BDS}\left(\rho_{1}^{n},\,\V v^{\star,n+\myhalf},\,\grad\cdot\V F^{\star,n+\myhalf},\,\D t\right) & \mbox{ for BDS}\\
\frac{1}{2}\left(\rho_{1}\V v\right)^{n}+\frac{1}{2}\left(\rho_{1}\V v\right)^{\star,n+1} & \mbox{ for centered}
\end{cases},
\end{eqnarray*}
where $\V F^{\star,n+\myhalf}=\left(\V F^{n}+\V F^{\star,n+1}\right)/2$
and $\V v^{\star,n+\myhalf}=\left(\V v^{n}+\V v^{\star,n+1}\right)/2$.
\item \label{enu:InertialCorrector}Compute $c^{n+1}=\rho_{1}^{n+1}/\rho^{n+1}$
and compute
\[
\V F^{n+1}=\left(\rho\chi\grad c\right)^{n+1}+\sqrt{\frac{2\left(\chi\rho\mu_{c}^{-1}\right)^{n+1}k_{B}T}{\D t\,\D V}}\,\widetilde{\V W}^{n+1}.
\]

\item Take a corrector step for velocity by solving the Stokes system, using
$\V v^{\star,n+1}$ as a reference state,
\begin{eqnarray*}
\frac{\rho^{n+1}\V v^{n+1}-\rho^{n}\V v^{n}}{\D t}+\grad\pi^{n+\frac{1}{2}} & = & \frac{1}{2}\grad\cdot\left(\left(-\rho\V v\V v\right)^{n}+\left(-\rho\V v\V v\right)^{\star,n+1}\right)+\frac{1}{2}\left(\rho^{n}+\rho^{n+1}\right)\V g\\
 & + & \frac{1}{2}\grad\cdot\left(\eta^{n}\bar{\grad}\V v^{n}+\eta^{n+1}\bar{\grad}\V v^{n+1}\right)+\frac{1}{2}\grad\cdot\left[\left(\sqrt{\frac{\eta^{n}k_{B}T}{\D t\,\D V}}+\sqrt{\frac{\eta^{n+1}k_{B}T}{\D t\,\D V}}\right)\,\overline{\M W}^{n}\right]\\
\grad\cdot\V v^{n+1} & = & -\grad\cdot\left(\beta\rho^{-1}\V F^{n+1}\right).
\end{eqnarray*}
.\end{enumerate}
\end{algorithm}

Several variants of the inertial Algorithm \ref{alg:LMInertial} preserve
deterministic second-order accuracy. For example, in the corrector
stage for $\rho_{1}$, for centered advection we use a trapezoidal
approximation to the advective flux,
\begin{equation}
\frac{1}{2}\left(\rho_{1}\V v\right)^{n}+\frac{1}{2}\left(\rho_{1}\V v\right)^{\star,n+1},\label{eq:adv_trap}
\end{equation}
but we could have also used a midpoint approximation
\begin{equation}
\left(\frac{\rho_{1}^{n}+\rho_{1}^{\star,n+1}}{2}\right)\left(\frac{\V v^{n}+\V v^{\star,n+1}}{2}\right).\label{eq:adv_mid}
\end{equation}
without affecting the second-order weak accuracy \cite{MultiscaleIntegrators}.
Note that BDS advection by construction requires a midpoint approximation
to the advective velocity; no analysis of the order of stochastic
accuracy is available for BDS advection at present. In the corrector
step for velocity, in Algorithm \ref{alg:LMInertial} we use corrected
values for the viscosity, but one can also use the values from the
predictor $\eta^{*,n+1}$.

\begin{algorithm}
\caption{\label{alg:LMOverdamped}A time step of our implicit midpoint temporal
integrator for the overdamped equations (\ref{eq:rho1_eq},\ref{eq:rho_eq},\ref{eq:steady_Stokes}).}

\begin{enumerate}
\item Calculate predictor diffusive fluxes and generate stochastic fluxes
for a half step to the midpoint,
\[
\V F^{n}=\left(\rho\chi\grad c\right)^{n}+\sqrt{\frac{2\left(\chi\rho\mu_{c}^{-1}\right)^{n}k_{B}T}{\left(\D t/2\right)\,\D V}}\,\widetilde{\V W}_{A}^{n}.
\]

\item Generate a random advection velocity by solving the steady Stokes
equation with random forcing,
\begin{eqnarray*}
\grad\pi^{n} & = & \grad\cdot\left(\eta^{n}\bar{\grad}\V v^{n}\right)+\grad\cdot\left(\sqrt{\frac{\eta^{n}k_{B}T}{\left(\D t/2\right)\,\D V}}\,\overline{\M W}_{A}^{n}\right)+\rho^{n}\V g\\
\grad\cdot\V v^{n} & = & -\grad\cdot\left(\beta\rho^{-1}\V F^{n}\right).
\end{eqnarray*}

\item Take a predictor midpoint Euler step for $\rho_{1}$, and similarly
for $\rho_{2}$, or, equivalently, for $\rho$,
\begin{eqnarray*}
\rho_{1}^{\star,n+\myhalf} & = & \rho_{1}^{n}+\frac{\D t}{2}\,\grad\cdot\V F^{n}-\frac{\D t}{2}\,\grad\cdot\begin{cases}
\text{BDS}\left(\rho_{1}^{n},\,\V v^{n},\,\grad\cdot\V F^{n},\,\frac{\D t}{2}\right) & \mbox{ for BDS}\\
\rho_{1}^{n}\V v^{n} & \mbox{ for centered}
\end{cases},
\end{eqnarray*}
and compute $c^{\star,n+\myhalf}=\rho_{1}^{\star,n+\myhalf}/\rho^{\star,n+\myhalf}$.
\item Calculate corrector diffusive fluxes and generate stochastic fluxes,
\[
\V F^{\star,n+\myhalf}=\left(\rho\chi\grad c\right)^{\star,n+\myhalf}+\sqrt{\frac{2\left(\chi\rho\mu_{c}^{-1}\right)^{\star,n+\myhalf}k_{B}T}{\D t\,\D V}}\,\left(\frac{\widetilde{\V W}_{A}^{n}+\widetilde{\V W}_{B}^{n}}{\sqrt{2}}\right),
\]
where $\widetilde{\V W}_{B}^{n}$ is a collection of random numbers
generated independently of $\widetilde{\V W}_{A}^{n}$.
\item Solve the corrected steady Stokes equation
\begin{eqnarray*}
\grad\pi^{\star,n+\myhalf} & = & \grad\cdot\left(\eta^{\star,n+\myhalf}\bar{\grad}\V v^{\star,n+\myhalf}\right)+\grad\cdot\left[\sqrt{\frac{\eta^{\star,n+\myhalf}k_{B}T}{\D t\,\D V}}\,\left(\frac{\overline{\M W}_{A}^{n}+\overline{\M W}_{B}^{n}}{\sqrt{2}}\right)\right]+\rho^{\star,n+\myhalf}\V g\\
\grad\cdot\V v^{\star,n+\myhalf} & = & -\grad\cdot\left(\beta\rho^{-1}\V F^{\star,n+\myhalf}\right).
\end{eqnarray*}

\item Correct $\rho_{1}$, and similarly for $\rho_{2}$, or, equivalently,
for $\rho$, 
\begin{eqnarray*}
\rho_{1}^{n+1} & = & \rho_{1}^{n}+\D t\,\grad\cdot\V F^{\star,n+\myhalf}-\D t\,\grad\cdot\begin{cases}
\text{BDS}\left(\rho_{1}^{n},\,\V v^{\star,n+\myhalf},\,\grad\cdot\V F^{\star,n+\myhalf},\,\D t\right) & \mbox{ for BDS}\\
\left(\rho\V v\right)^{\star,n+\myhalf} & \mbox{ for centered}
\end{cases},
\end{eqnarray*}
and set $c^{n+1}=\rho_{1}^{n+1}/\rho^{n+1}$.\end{enumerate}
\end{algorithm}

\subsubsection{Order of Accuracy}

For explicit temporal integrators, we relied on a gauge formulation
to write the low Mach equations in the form of a standard unconstrained
initial-value problem, thus allowing us to use standard integrators
for ODEs \cite{LowMachExplicit}. In the semi-implicit case, however,
we do not use a gauge formulation because the Stokes solver we use
works directly with the pressure and velocity. This makes proving
second-order temporal accuracy nontrivial even in the deterministic
context; we therefore rely on empirical convergence testing to confirm
the second-order deterministic accuracy.

In the stochastic context, there is presently no available theoretical
analysis when BDS advection is employed; existing analysis \cite{LLNS_S_k,DFDB,MultiscaleIntegrators}
assumes a method of lines (MOL) discretization in which space is discretized
first to obtain a system of SODEs. For centered advection, which does
lead to an MOL discretization, the algorithms used here are based
on the second-order weak temporal integrators developed in Ref. \cite{MultiscaleIntegrators}.
In particular, for the case of the inertial equations (\ref{eq:momentum_eq},\ref{eq:rho1_eq},\ref{eq:div_v_constraint},\ref{eq:rho_eq})
we base our temporal integrator on an implicit trapezoidal method.
It should be emphasized however that the analysis in Ref. \cite{MultiscaleIntegrators}
applies to unconstrained Langevin systems, while the low Mach equations
are constrained by the EOS. Nevertheless, the deterministic accuracy
of the method is crucial even when fluctuations of primary interest
because in linearized fluctuating hydrodynamics the fluctuations are
linearized around the solution of the determinisitc equations, which
must itself be computed numerically \cite{MultiscaleIntegrators}
accurately in order to have any chance to compute the fluctuations
accurately. For the case of the overdamped equations (\ref{eq:rho1_eq},\ref{eq:rho_eq},\ref{eq:steady_Stokes}),
we base our temporal integrator on an implicit midpoint method. In
this case the analysis presented in Ref. \cite{MultiscaleIntegrators}
does apply since the velocity is not a variable in the overdamped
equations and the limiting equation for concentration is unconstrained.
This analysis indicates that the overdamped temporal integrator below
is second-order weakly accurate for the \emph{linearized} overdamped
low Mach number equations.

\section{\label{sub:Validation}Validation and Testing}

In this section we apply the inertial and overdamped low Mach algorithms
described in Section \ref{sec:SpatioTemporal} in a stochastic and
several deterministic contexts. First, we demonstrate our ability
to accurately model equilibrium fluctuations by analyzing the static
spectrum of the fluctuations. Next, we confirm the second-order deterministic
order of accuracy of our methods on a low Mach number lid-driven cavity
test. Next, we confirm that the BDS advection scheme enables robust
simulation in cases when there is little or no mass diffusion (i.e.,
nearly infinite Peclet number). Lastly, we use the inertial algorithm
to simulate the development of a Kelvin-Helmholtz instability when
a lighter less viscous fluid is impulsively set in motion on top of
a heavier more viscous fluid.

\subsection{Equilibrium Fluctuations}

One of the key quantities used to characterize the intensity of \textit{equilibrium}
thermal fluctuations is the static structure factor or static spectrum
of the fluctuations at thermodynamic equilibrium. We examine the static
structure factors in both the inertial and overdamped regimes. We
use arbitrary units with $T=1$, $k_{B}=1$, molecular masses $m_{1}=1,$
$m_{2}=2$, and pure component densities $\bar{\rho}_{1}=2/3$, $\bar{\rho}_{2}=2$.
We initialize the domain with $c=0.5,$ which gives $\rho=1$. The
diffusion coefficient was constant$\chi=1$, whereas the viscosity
varies linearly from $\eta=1$ to $\eta=10$ (for the inertial tests),
and from $\eta=1$ to $\eta=100$ (for the overdamped tests) as $c$
varies from $0$ to $1$, but note that at equilibrium the concentration
fluctuations are small so the viscosity varies little over the domain.
We assume an ideal mixture, giving chemical potential $\mu_{c}^{-1}k_{B}T=c(1-c)[cm_{2}+(1-c)m_{1}]$
\cite{LowMachExplicit,Bell:09}. At these conditions, the equilibrium
density variance is $\Delta V\langle(\delta\rho)^{2}\rangle=S_{\rho}=0.375$,
where $\Delta V$ is the volume of a grid cell (see Appendix A1 in
\cite{LowMachExplicit}). We use a periodic system with $32\times32$
grid cells with $\Delta x=\Delta y=1$, with the thickness in the
third direction set to give a large $\Delta V=10^{6}$ and thus small
fluctuations, ensuring consistency with linearized fluctuating hydrodynamics.
A total of $10^{5}$ time steps are skipped in the beginning to allow
equilibration of the system, and statistics are then collected for
an additional $10^{6}$ steps. We run both the inertial and overdamped
algorithms using three different time steps, $\Delta t=0.1,$ $0.05,$
and $0.025$, the largest of which corresponds to 40\% of the maximum
allowable time step by the explicit mass diffusion CFL condition.

\begin{table}[h]
\begin{centering}
\caption{\label{tab:Srho} Equilibrium static structure factor $S_{\rho}$
averaged over all wavevectors for the inertial and overdamped algorithms
using three different time steps. The exact solution from theory is
$S_{\rho}=0.375$, allowing us to estimate an order of accuracy from
the average error over all wavenumbers. Note that there are significant
statistical errors present, especially at small wavenumbers, and these
make it difficult to reliably estimate the asymptotic order of accuracy
empirically when the error is very small (as for the overdamped integrator).}

\par\end{centering}

\centering{}%
\begin{tabular}{ccccc}
 &
$\Delta t$ &
$S_{\rho}$ &
|Error|  &
Order\tabularnewline
\hline 
Inertial &
0.1 &
0.3201 &
0.0549 &
\tabularnewline
 &
0.05 &
0.3624 &
0.0126 &
2.12\tabularnewline
 &
0.025 &
0.3722 &
0.0029 &
2.14\tabularnewline
\hline 
Overdamped  &
0.1 &
0.4192 &
0.0442 &
\tabularnewline
 &
0.05 &
0.3786 &
0.0036 &
3.63\tabularnewline
 &
0.025 &
0.3755 &
0.0005 &
2.92\tabularnewline
\end{tabular}
\end{table}

\begin{figure}
\centering{}\includegraphics[width=0.33\textwidth]{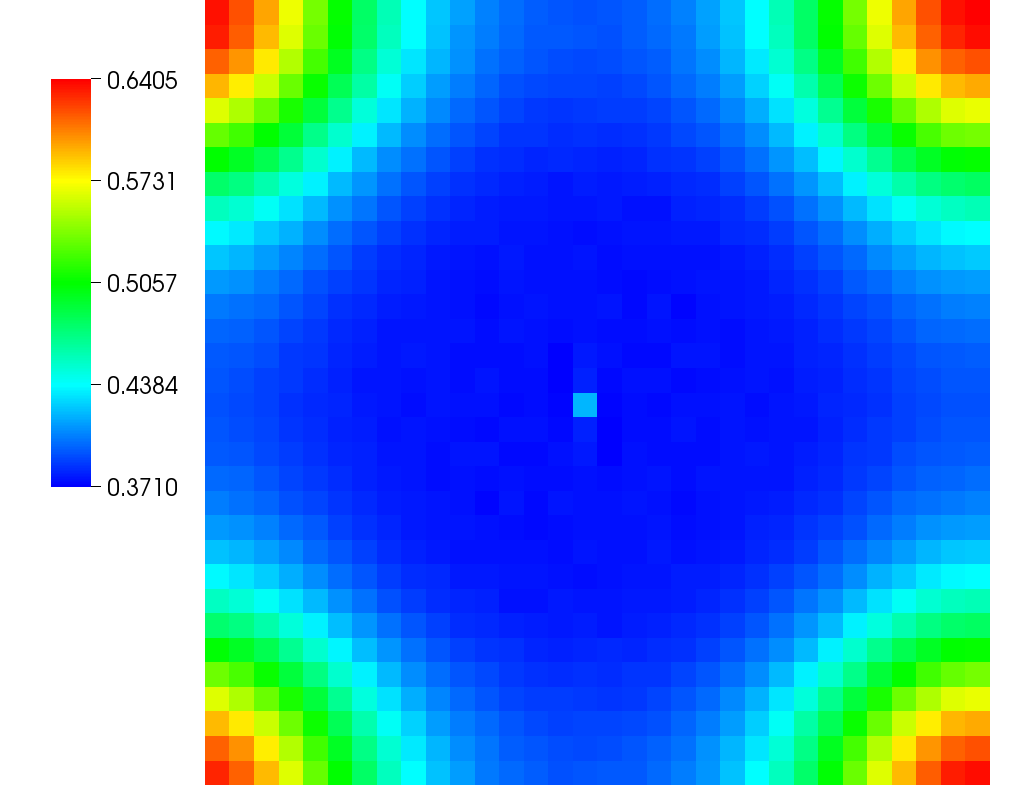}
\includegraphics[width=0.33\textwidth]{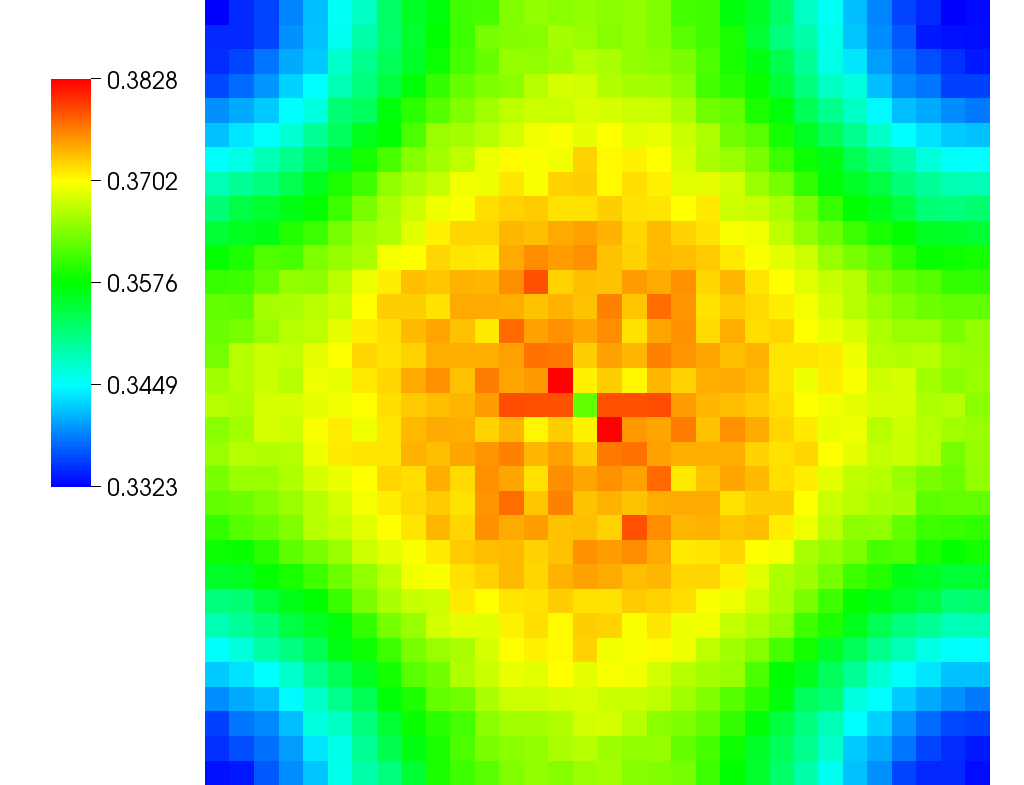}\includegraphics[width=0.33\textwidth]{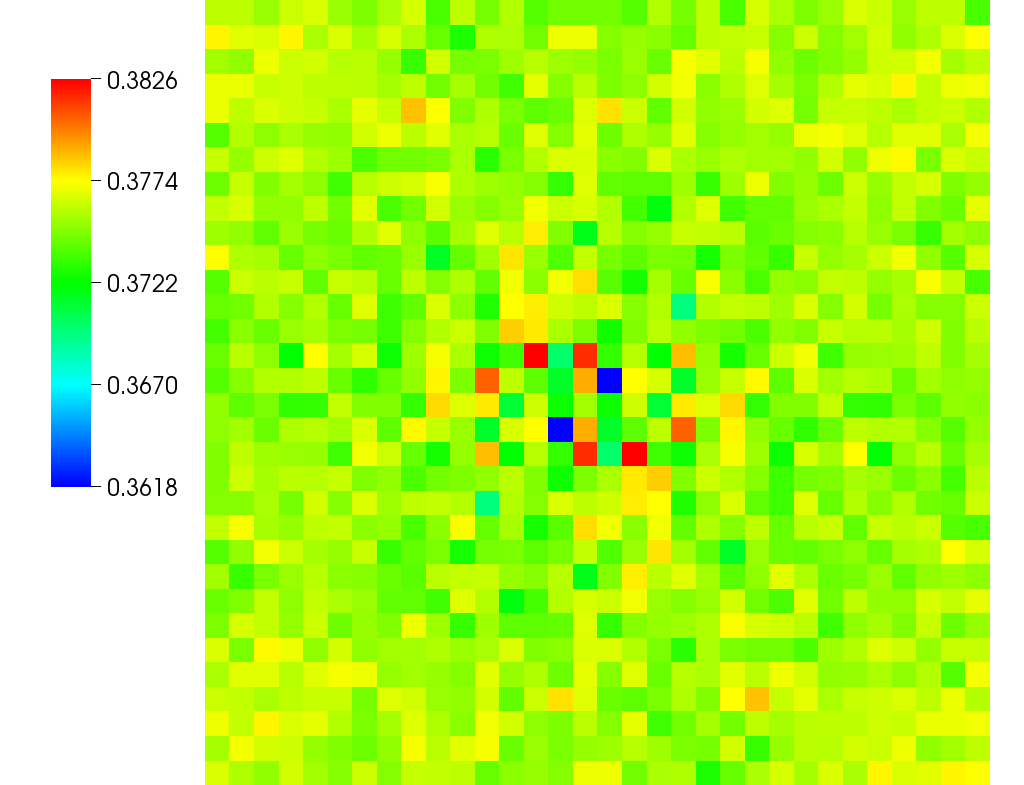}
\caption{\label{fig:Sk_overdamped}Equilibrium static structure factor $S_{\rho}$
as a function of wavevector (zero being at the center of the figures)
for the overdamped simulations with $\Delta t=0.1$ (left), $\Delta t=0.05$
(middle), and $\Delta t=0.025$ (right). The correct result, which
is recovered in the limit $\Delta t\rightarrow0$, is $S_{\rho}=0.375$.
The artifacts decrease by roughly a factor of 8 as the time step is
reduced in half.}
\end{figure}

In Table \ref{tab:Srho} we observe that as we reduce the time step
by a factor of two, we see a reduction in error in the average value
of $S_{\rho}$ over all wavenumbers by a factor of \textasciitilde{}4
(second-order convergence) for the inertial algorithm, and a factor
of \textasciitilde{}8 (third-order convergence) for the overdamped
algorithm (the latter being consistent with the fact that the explicit
midpoint method is third-order accurate for static covariances \cite{DFDB}).
In Figure \ref{fig:Sk_overdamped} we show the spectrum of density
fluctuations at equilibrium for three different time step sizes. At
thermodynamic equilibrium, the static structure factors are independent
of the wavenumber due to the local nature of the correlations. Since
we include mass diffusion using an explicit temporal integrator, for
larger time steps we expect to see additional deviation from a flat
spectrum at the largest wavenumbers (i.e., for $k\sim\Delta x^{-1}$)\cite{LLNS_S_k,DFDB}.
In the limit of sufficiently small time steps, we recover the correct
flat spectrum, demonstrating that our model and numerical scheme obey
a fluctuation-dissipation principle.

\subsection{Deterministic Lid-Driven Cavity Convergence Test}

In this section, we simulate a smooth test problem and empirically
confirm deterministic second-order accuracy of Algorithms \ref{alg:LMInertial}
and \ref{alg:LMOverdamped} even in the presence of boundary conditions,
inertial effects, and gravity, as well as nonconstant density, mass
diffusion coefficient, and viscosity. The problem is a deterministic
lid-driven cavity flow, following previous work by Boyce Griffith
for incompressible constant-density and constant-viscosity flow \cite{NonProjection_Griffith}.

\begin{figure}
\centering{}\includegraphics[width=2.9in]{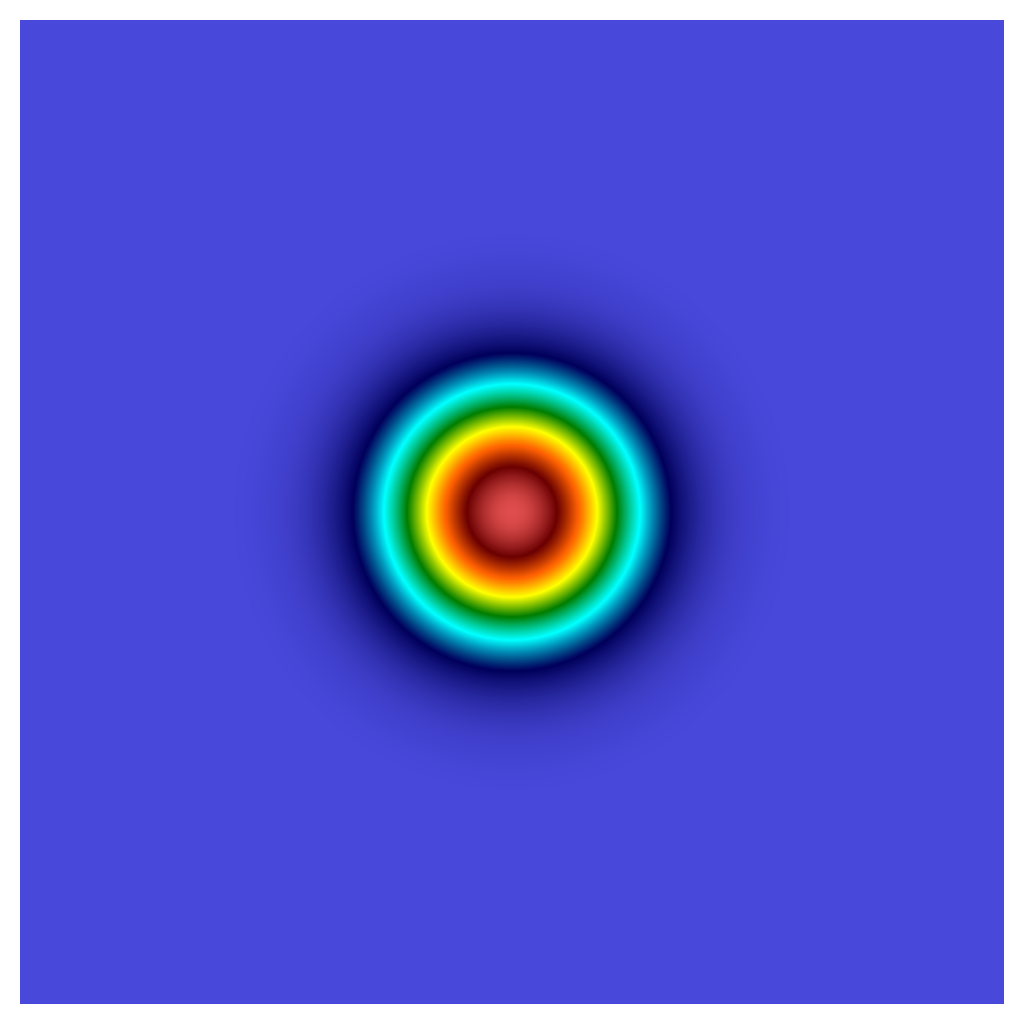}
\includegraphics[width=2.9in]{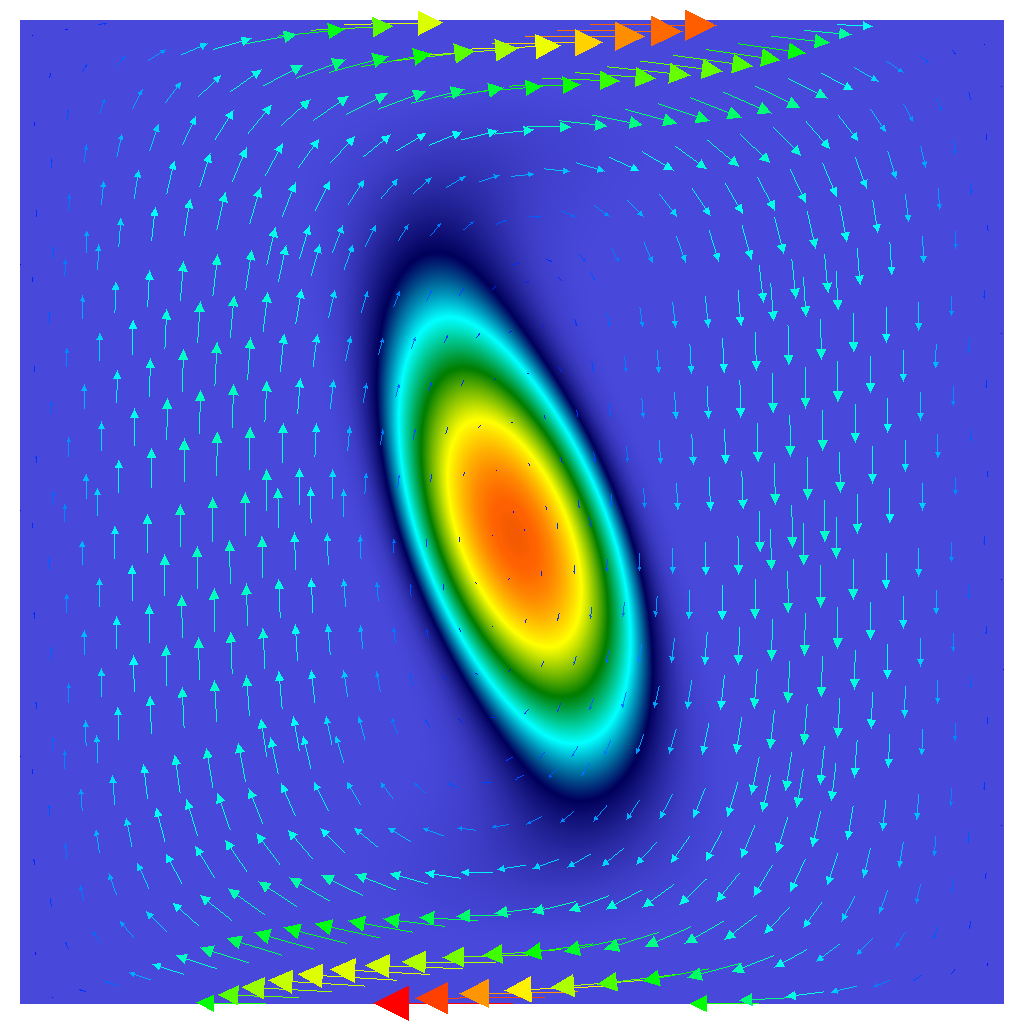}
\caption{\label{fig:LidDriven}Initial ($t=0$) and final ($t=2$) concentration
(scalar color field) and velocities (vector field) for the low Mach
number lid-driven cavity test problem.}
\end{figure}

We use CGS units (centimeters for length, seconds for time, grams
for mass). We consider a square (two dimensions) or cubic (three dimensions)
domain with side of length $L=1$ bounded on all sides by no-slip
walls moving with a specified velocity. The bottom and top walls ($y$-direction)
are no-slip walls moving in equal and opposite directions, setting
up a circular flow pattern, while the remaining walls are stationary.
The top wall has a specified velocity, in two dimensions, 
\begin{equation}
u(x,t)=\begin{cases}
\frac{1}{4}\left[1+\sin{(2\pi x-\frac{\pi}{2}})\right]\left[1+\sin{(2\pi t-\frac{\pi}{2}})\right], & t<1/2\\
\frac{1}{2}\left[1+\sin{(2\pi x-\frac{\pi}{2}})\right], & t\ge1/2
\end{cases}
\end{equation}
and in three dimensions,
\begin{equation}
u(x,z,t)=w(x,z,t)=\begin{cases}
\frac{1}{8}\left[1+\sin{(2\pi x-\frac{\pi}{2}})\right]\left[1+\sin{(2\pi z-\frac{\pi}{2}})\right]\left[1+\sin{(2\pi t-\frac{\pi}{2}})\right], & t<1/2\\
\frac{1}{4}\left[1+\sin{(2\pi x-\frac{\pi}{2}})\right]\left[1+\sin{(2\pi z-\frac{\pi}{2}})\right]. & t\ge1/2
\end{cases}.
\end{equation}
Note that the wall velocity tapers to zero at the corners in order
to regularize the corner singularities \cite{NonProjection_Griffith};
similarly, the velocity smoothly increases with time to its final
value in order to avoid potential loss of accuracy due to an impulsive
start of the flow. The two liquids have pure-component densities $\bar{\rho}_{1}=2$
and $\bar{\rho}_{2}=1$. The initial conditions are $\vb=0$ for velocity,
and a Gaussian bump of higher density for the concentration, $c\left(\V r,t\right)=\exp\left(-75r^{2}\right)$,
where $r$ is the distance to the center of the domain. The viscosity
varies linearly as a function of concentration, such that $\eta=0.1$
when $c=0$ and $\eta=1$ when $c=1$. Similarly, the mass diffusion
coefficient varies linearly as a function of concentration, such that
$\chi=10^{-4}$ when $c=0$ and $\chi=10^{-3}$ when $c=1$. In order
to confirm that second-order accuracy is preserved even in the limit
of infinite Peclet number if BDS advection is employed, we also perform
simulations with $\chi=0$. Figure \ref{fig:LidDriven} illustrates
the initial and final (at time $t=2$) configurations of concentration
and velocity in two dimensions.

Recall that advection of the concentration can be treated using centered
advection or the BDS advection scheme (see Section \ref{sub:BDS}).
BDS advection can use either a bi-linear (trilinear in 3D) or a quadratic
reconstruction (2D only), and can further be limited to avoid the
appearance of spurious local extrema. Here we present convergence
results for the following test problems: 
\begin{itemize}
\item Test 1: Centered advection, nonzero $\chi$
\item Test 2: Unlimited bilinear BDS advection, nonzero $\chi$
\item Test 3: Unlimited quadratic BDS advection, nonzero $\chi$
\item Test 4: Unlimited bilinear BDS advection, $\chi=0$.
\end{itemize}
We perform Tests 1-4 using both the inertial Algorithm \ref{alg:LMInertial}
and the overdamped Algorithm \ref{alg:LMOverdamped}. The Reynolds
number in this test is of order unity and there is only a small difference
in the results for the inertial and overdamped equations. Recall that
concentration is the \emph{only} independent variable in the overdamped
equations.

In two dimensions, we discretize the problem on a grid of $64^{2}$,
$128^{2}$, $256^{2}$ or $512^{2}$ grid cells. The time step size
for the coarsest simulation is $\Delta t=5\times10^{-3}$ and is reduced
by a factor of 2 as the resolution increases by a factor of 2. This
corresponds to an advective Courant number of $v_{\text{max}}\D t/\D x\sim0.3$
for each simulation. The diffusive Courant number is $\chi\D t/\D x^{2}\sim0.16$
(recall that the stability limit is $1/4=0.25$ in two dimensions)
for the finest simulation, reducing by a factor of 2 with each successive
grid coarsening. We simulate the flow and compute error norms at time
$t=2$. In table \ref{tab:Linf-inertial} we present estimates of
the order of convergence in the $L_{\infty}$ (max) norm for the velocity
components and concentration for the inertial equations. We see clear
second-order pointwise convergence, without any artifacts near the
boundaries. Similar results are obtained for the concentration in
the overdamped limit, as shown in table \ref{tab:Linf-overdamped}.

\begin{table}[h]
\begin{centering}
\caption{\label{tab:Linf-inertial}Convergence of errors in the $L_{\infty}$
norm for a \emph{two-dimensional} \emph{inertial} low Mach lid-driven
cavity problem as the grid is refined in space and time, for the components
of the velocity $\V v=\left(u,v\right)$ and concentration $c$. Order
of convergence is estimated from the error ratio between two successive
refinements.}

\par\end{centering}

\centering{}%
\begin{tabular}{ccccccc}
refinement &
 &
64-128  &
Order  &
128-256  &
Order  &
256-512 \tabularnewline
\hline 
Test 1:  &
$u$  &
1.93e-03  &
1.91  &
5.12e-04  &
1.96  &
1.32e-04 \tabularnewline
 &
$v$  &
8.69e-04  &
1.99  &
2.19e-04  &
2.00  &
5.49e-05 \tabularnewline
 &
$c$  &
3.02e-04  &
1.99  &
7.60e-04  &
2.00  &
1.90e-04 \tabularnewline
\hline 
Test 2:  &
$u$  &
1.92e-03  &
1.91  &
5.11e-04  &
1.95  &
1.32e-04 \tabularnewline
 &
$v$  &
9.08e-04  &
1.96  &
2.34e-04  &
1.99  &
5.91e-05 \tabularnewline
 &
$c$  &
2.63e-03  &
1.72  &
7.99e-04  &
1.92  &
2.11e-04 \tabularnewline
\hline 
Test 3:  &
$u$  &
1.92e-03  &
1.91  &
5.11e-04  &
1.95  &
1.32e-04 \tabularnewline
 &
$v$  &
8.62e-04  &
1.95  &
2.23e-04  &
1.98  &
5.64e-05 \tabularnewline
 &
$c$  &
1.95e-03  &
1.99  &
4.91e-04  &
2.00  &
1.23e-04 \tabularnewline
\hline 
Test 4:  &
$u$  &
1.91e-03  &
1.92  &
5.06e-04  &
1.96  &
1.30e-04 \tabularnewline
 &
$v$  &
9.78e-04  &
2.01  &
2.43e-04  &
2.02  &
6.00e-05 \tabularnewline
 &
$c$  &
4.29e-03  &
1.90  &
1.15e-03  &
1.97  &
2.93e-04 \tabularnewline
\end{tabular}
\end{table}

\begin{table}[h]
\begin{centering}
\caption{\label{tab:Linf-overdamped}Convergence of errors in the $L_{\infty}$
norm for a \emph{two-dimensional} \emph{overdamped} low Mach lid-driven
cavity problem as the grid is refined in space and time, for concentration
$c$. Order of convergence is estimated from the error ratio between
two successive refinements.}

\par\end{centering}

\centering{}%
\begin{tabular}{cccccc}
refinement &
64-128  &
Order  &
128-256  &
Order  &
256-512 \tabularnewline
\hline 
Test 1:  &
3.57e-03  &
2.01  &
8.89e-04  &
2.00  &
2.22e-04 \tabularnewline
\hline 
Test 2:  &
2.70e-03  &
1.78  &
7.87e-04  &
1.92  &
2.08e-04 \tabularnewline
\hline 
Test 3:  &
1.95e-03  &
1.98  &
4.95e-04  &
1.89  &
1.34e-04 \tabularnewline
\hline 
Test 4:  &
4.23e-03  &
1.93  &
1.11e-03  &
1.96  &
2.86e-04 \tabularnewline
\end{tabular}
\end{table}

In three dimensions, we discretize the problem on a grid of $32^{3}$,
$64^{3}$, $128^{3}$ or $256^{3}$ grid cells. The time step size
for the coarsest simulation is $\Delta t=1.25\times10^{-2}$, which
corresponds to an advective Courant number of $\sim0.4$, and diffusive
Courant number of $\sim0.10$ (stability limit is $1/6\approx0.17$)
for the finest resolution simulation. We simulate the flow and compute
error norms at time $t=1$. We limit our study here to inertial flow
and only perform Tests 1 and 2 (note that there is presently no available
unlimited quadratic BDS advection scheme in three dimensions, so test
3 cannot be performed). We also try a higher-order one-sided difference
for the tangential velocity at the no-slip boundaries, which does
not affect the asymptotic rate of convergence, but it can significantly
reduce the magnitude of the errors, and enables us to reach the asymptotic
regime for smaller grid sizes \cite{NonProjection_Griffith}. The
numerical convergence results shown in table \ref{tab:Linf_3d} demonstrate
the second-order deterministic accuracy of our method in three dimensions.

\begin{table}[h]
\begin{centering}
\caption{\label{tab:Linf_3d}Convergence of errors in the $L_{\infty}$ norm
for a \emph{three-dimensional} \emph{inertial} low Mach lid-driven
cavity problem as the grid is refined in space and time, for the components
of the velocity $\V v=\left(u,v,w\right)$ and concentration $c$.
Order of convergence is estimated from the error ratio between two
successive refinements.}

\par\end{centering}

\centering{}%
\begin{tabular}{ccccccc}
 &
 &
32-64  &
Rate  &
64-128  &
Rate  &
128-256 \tabularnewline
\hline 
Test 1:  &
$u$  &
7.66e-03  &
1.75  &
2.27e-03  &
1.88  &
6.16e-04 \tabularnewline
 &
$v$  &
3.12e-03  &
1.96  &
8.02e-04  &
1.99  &
2.02e-04 \tabularnewline
 &
$w$  &
7.66e-03  &
1.75  &
2.27e-03  &
1.88  &
6.16e-04 \tabularnewline
 &
$c$  &
1.22e-02  &
2.00  &
3.06e-03  &
2.00  &
7.64e-04 \tabularnewline
\hline 
Test 1:  &
$u$  &
2.30e-03  &
1.97  &
5.88e-04  &
2.02  &
1.45e-04 \tabularnewline
with higher-order &
$v$  &
9.01e-04  &
2.23  &
1.92e-04  &
1.99  &
4.82e-05 \tabularnewline
boundary stencil &
$w$  &
2.30e-03  &
1.97  &
5.88e-04  &
2.02  &
1.45e-04 \tabularnewline
 &
$c$  &
1.21e-02  &
1.99  &
3.05e-03  &
2.00  &
7.62e-04 \tabularnewline
\hline 
Test 2:  &
$u$  &
7.67e-03  &
1.75  &
2.28e-03  &
1.89  &
6.16e-04 \tabularnewline
 &
$v$  &
3.11e-03  &
1.96  &
8.01e-04  &
1.99  &
2.02e-04 \tabularnewline
 &
$w$  &
7.67e-03  &
1.75  &
2.28e-03  &
1.89  &
6.16e-04 \tabularnewline
 &
$c$  &
9.79e-03  &
1.91  &
2.61e-03  &
1.97  &
6.68e-04 \tabularnewline
\hline 
Test 2:  &
$u$  &
2.30e-03  &
1.96  &
5.89e-04  &
2.01  &
1.46e-04 \tabularnewline
with higher-order &
$v$  &
8.90e-04  &
2.21  &
1.92e-04  &
1.99  &
4.82e-05 \tabularnewline
boundary stencil &
$w$  &
2.30e-03  &
1.96  &
5.89e-04  &
2.01  &
1.46e-04 \tabularnewline
 &
$c$  &
9.70e-03  &
1.91  &
2.59e-03  &
1.96  &
6.67e-04 \tabularnewline
\end{tabular}
\end{table}

\subsection{Deterministic Sharp Interface Limit}

\begin{figure}
\centering{}\includegraphics[width=1.9in]{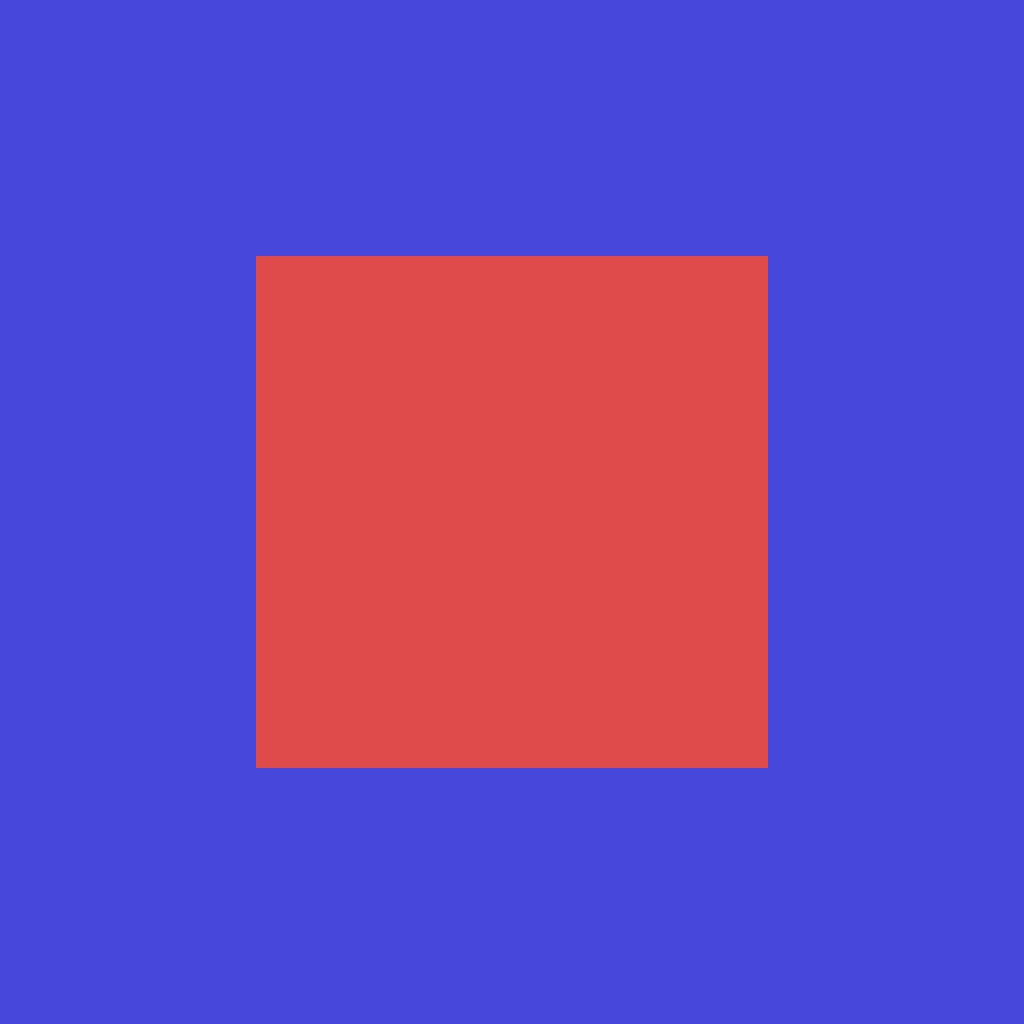}\hspace{1mm}\includegraphics[width=1.9in]{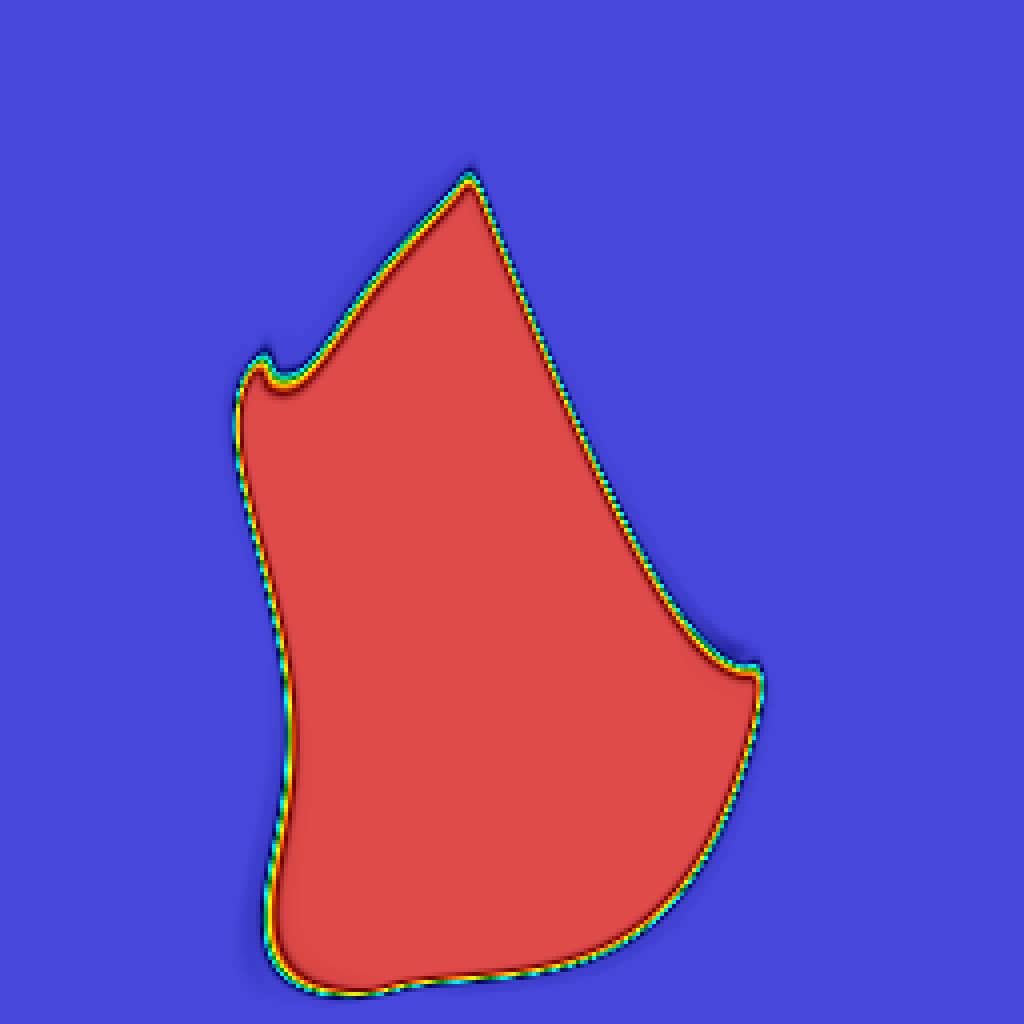}\hspace{1mm}\includegraphics[width=1.9in]{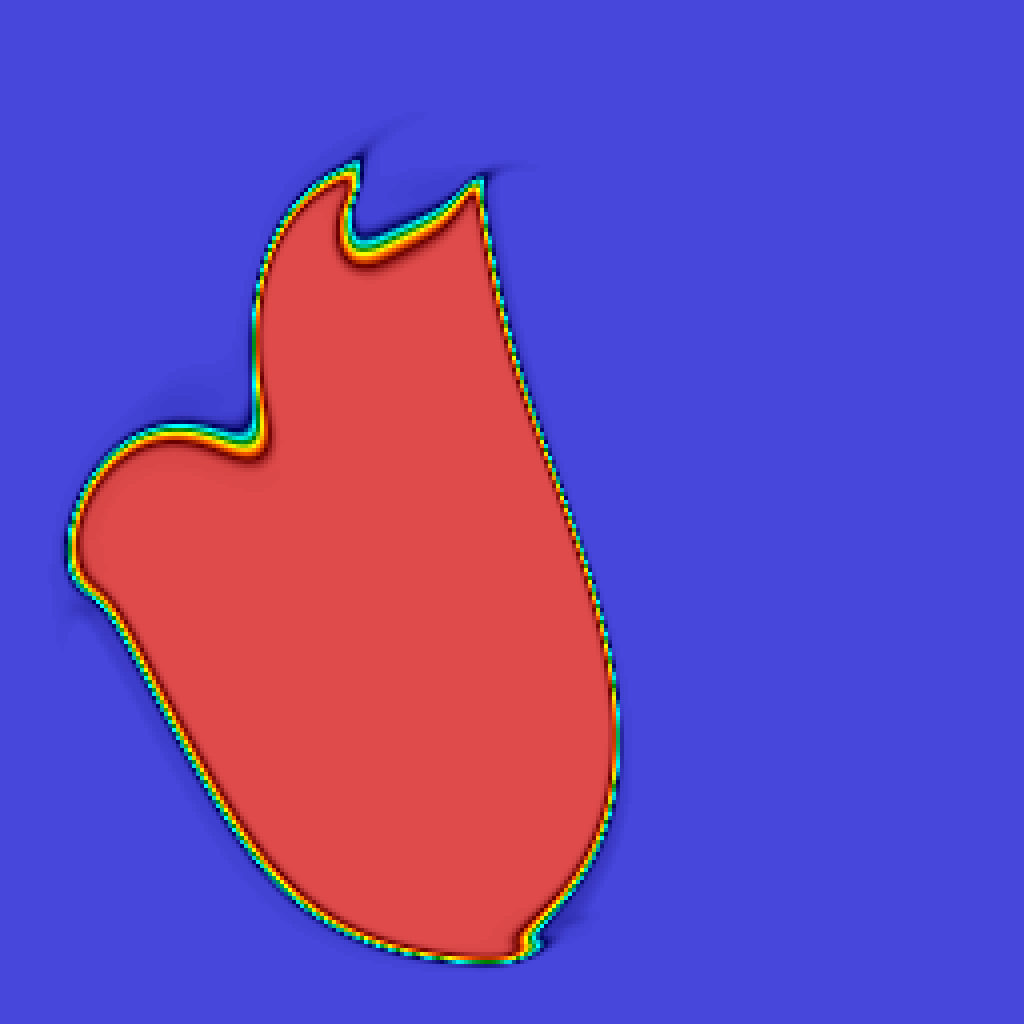}\hspace{1mm}\includegraphics[height=1.9in]{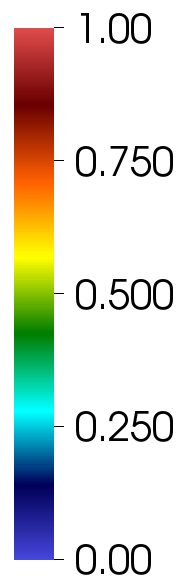}\caption{\label{fig:square}Advection of a square bubble in a lid-driven cavity
flow, using the limited quadratic BDS scheme. Concentration is shown
as a color plot at times $t=0,2,4$.}
\end{figure}

In this section we verify the ability of the BDS advection scheme
to advect concentration and density without creating spurious oscillations
or instabilities, even in the absence of mass diffusion, $\chi=0$,
and in the presence of sharp interfaces. The problem setup is similar
to the inertial lid-driven cavity test presented above, with the following
differences. First, the gravity is larger, $\V g=(0,-5)$, so that
the higher density region falls downward a significant distance. Secondly,
the initial conditions are a constant background of $c=0$ with a
square region covering the central 25\% of the domain initialized
to $c=1$ (see the left panel of Figure \ref{fig:square}). The correct
solution of the equations must remain a binary field, $c=1$ inside
the advected square curve, and $c=0$ elsewhere. In this test we employ
limited quadratic BDS, and use a grid of $256^{2}$ grid cells and
a fixed time step size $\Delta t=2.5\times10^{-3}$, corresponding
to an advective CFL number of $\sim0.6$. In Figure \ref{fig:square},
we show the concentration at several points in time, observing very
little smearing of the interface, even as the deformed bubble passes
near the bottom boundary.

\subsection{Deterministic Kelvin-Helmoltz Instability}

We simulate the development of a Kelvin-Helmholtz instability in three
dimensions in order to demonstrate the robustness of our inertial
time-advancement scheme in a deterministic setting. Our simulation
uses $256\times128\times256$ computational cells with grid spacing
$\Delta x=1/256$. We use periodic boundary conditions in the $x$
and $z$ directions, a no-slip conditions on the $y$ boundaries,
with prescribed velocity $\V v\left(x,y=0,z\right)=0$ on the bottom
boundary and $\V v\left(x,y=0.5,z\right)=(1,0,0)$ on the top boundary.
We use an adaptive time step size $\D t$ adjusted to maintain a maximum
advective CFL number $v_{max}\D t/\D x\leq0.9$. The binary fluid
mixture has a 10:1 density contrast with $\bar{\rho}_{1}=10$ and
$\bar{\rho}_{2}=1$. Viscosity varies linearly with concentration,
such that $\eta=10^{-4}$ for $c=0$ and $\eta=10^{-3}$ for $c=1$.
The mass diffusion coefficient is fixed at $\chi=10^{-6}$, which
makes the diffusive CFL number $\chi\D t/\D x^{2}\sim10^{-4}$, making
it necessary to use BDS advection in order to avoid instabilities
due to sharp gradients at the interface. We employ the bilinear BDS
advection scheme \cite{BDS} with limiting in order to preserve strict
monotonicity and maintain concentration within the bounds $0\leq c\leq1$.

\begin{figure*}
\begin{centering}
\includegraphics[width=0.49\textwidth]{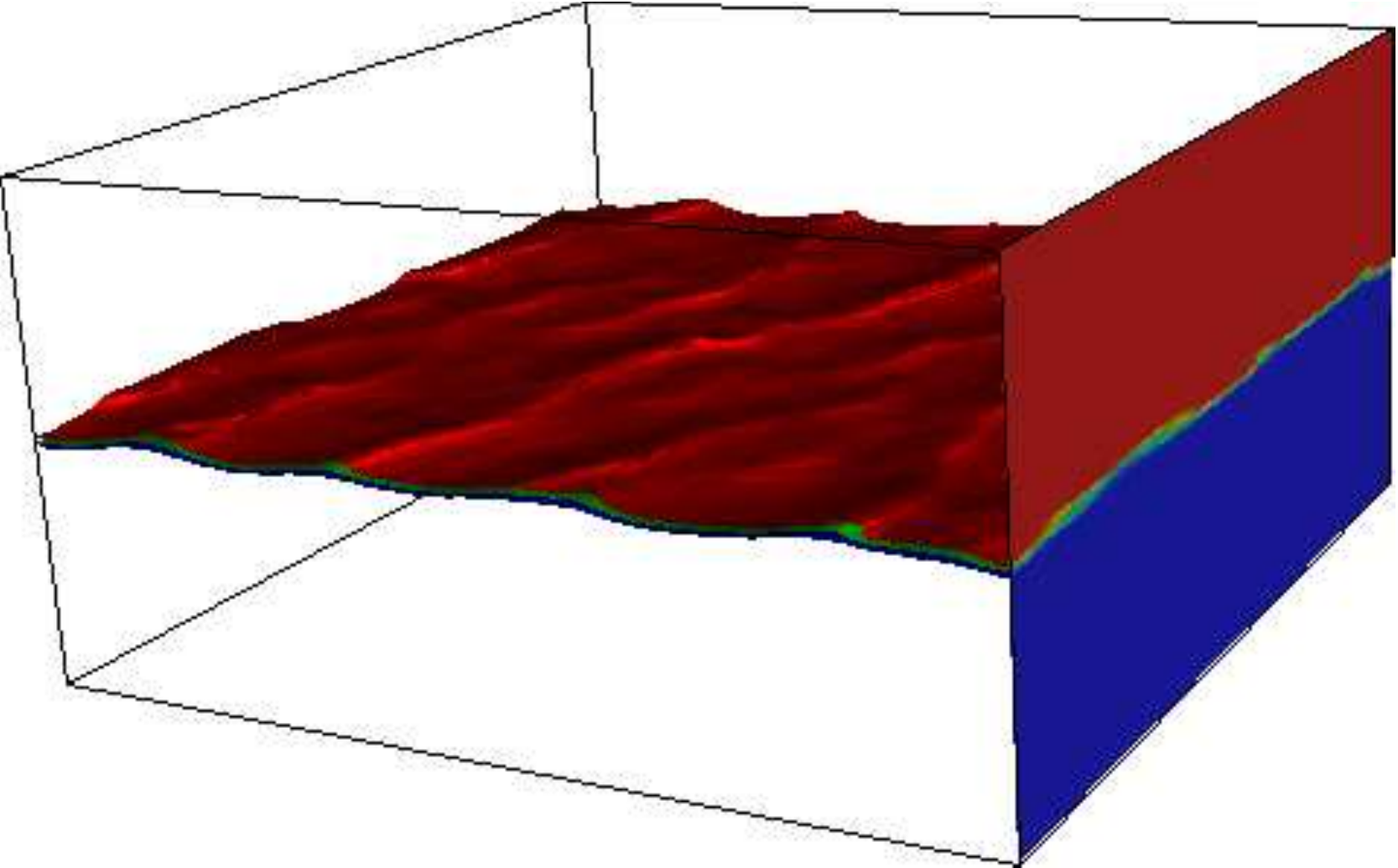}
\includegraphics[width=0.49\textwidth]{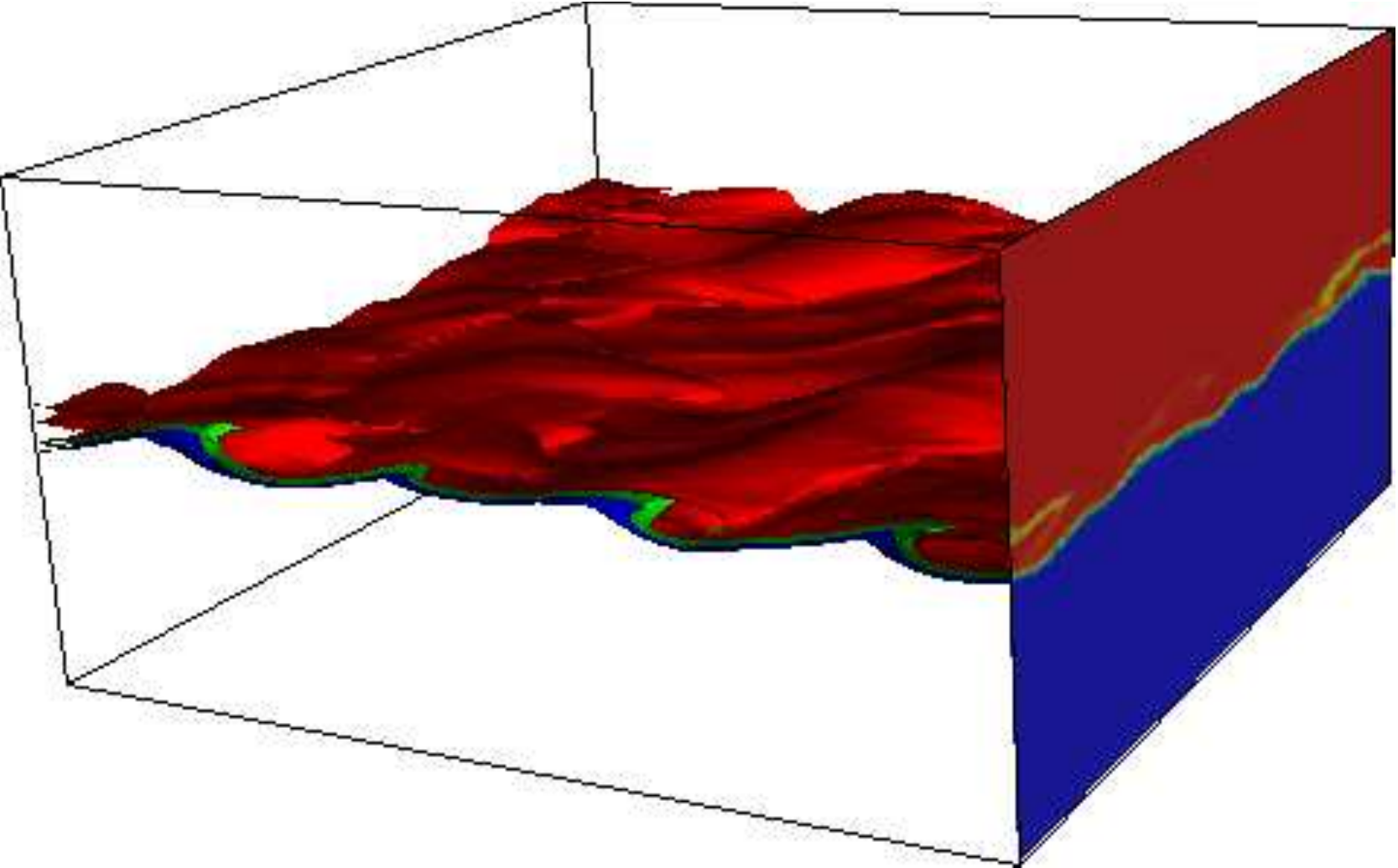}
\par\end{centering}

\centering{}\includegraphics[width=0.49\textwidth]{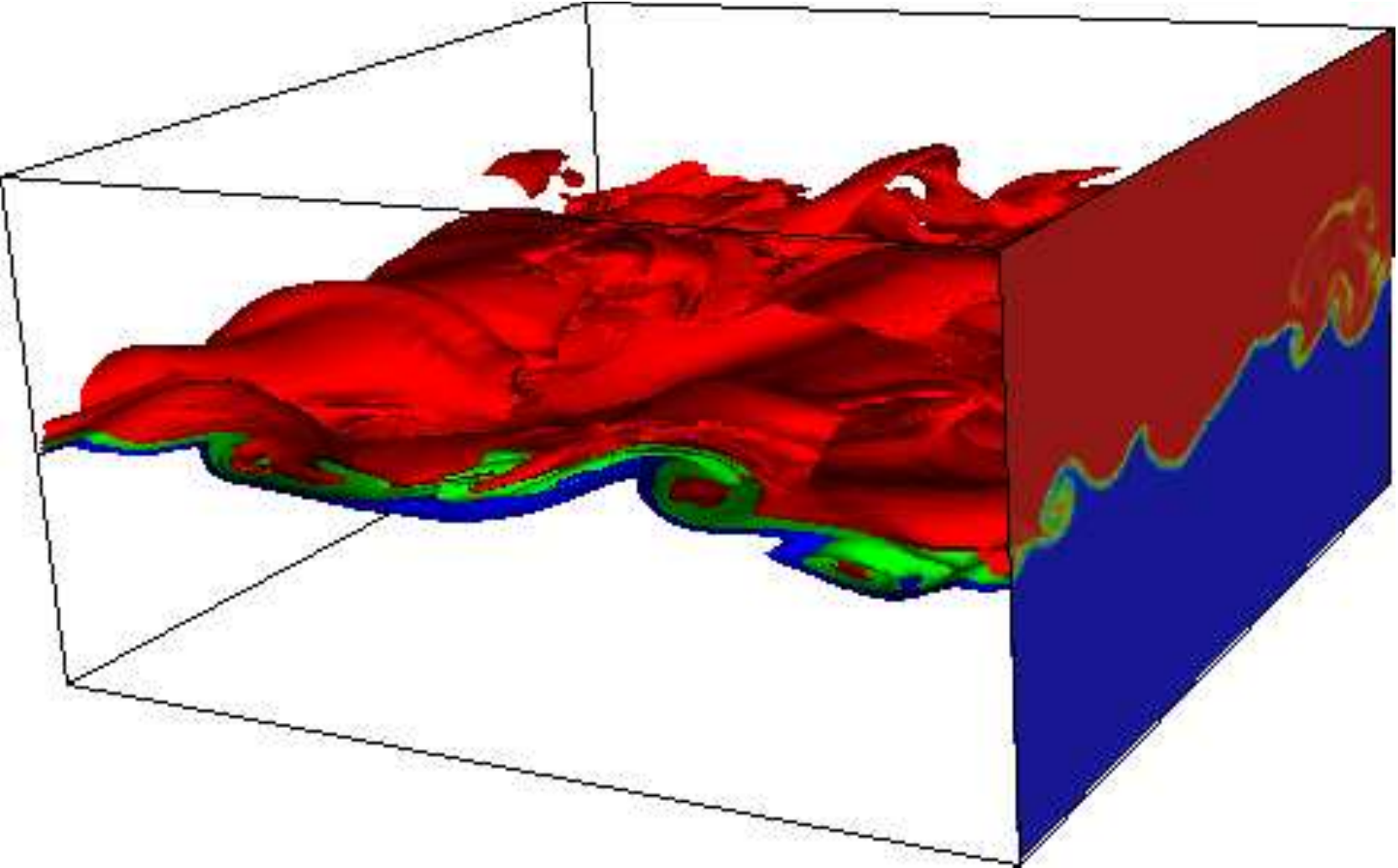}\includegraphics[width=0.49\textwidth]{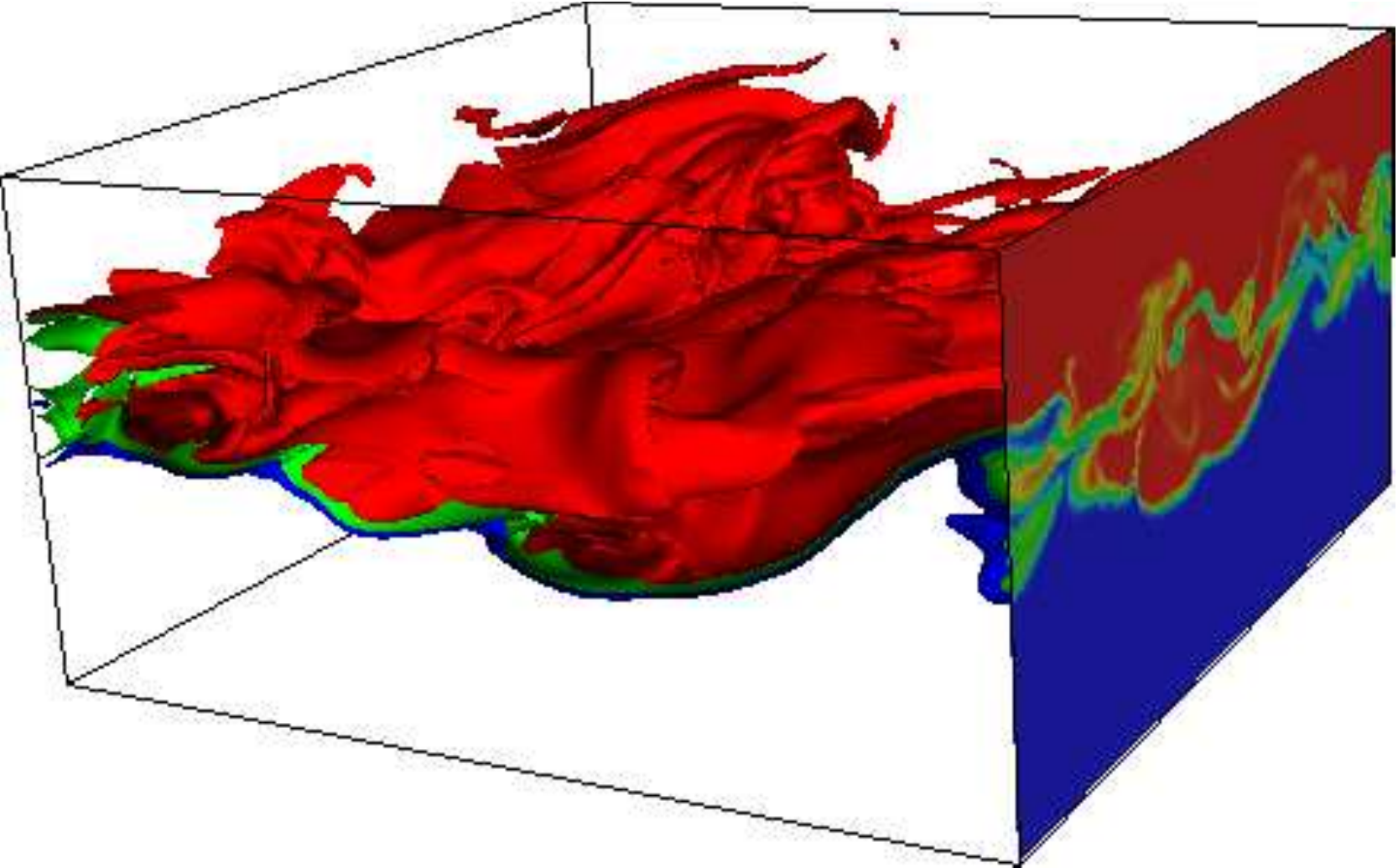}
\caption{\label{fig:KHLM_3D}The development of a Kelvin-Helmholtz instability
as a lighter less-viscous fluid streams over a ten times denser and
more viscous fluid. Contour surfaces of the density, ranging from
$\rho=1$ (red) to $\rho=10$ (blue), are shown at times $t=1.72,\,3.16,\,4.53$,
and $5.85$s.}
\end{figure*}

The initial condition is $c=1$ in the lower-half of the domain, and
$c=0$ in the upper-half of the domain, so that light fluid sits on
top of heavy fluid with a discontinuity in the concentration and velocity
at the interface. The initial momentum is set to $\rho\V v=(1,0,0)$
in the upper-half of the domain and $\rho\V v=0$ in the lower-half
of the domain. Gravity has a magnitude of $g=0.1$ acting in the downward
$y$-direction. In order to set off the instability, in a row of cells
at the centerline, $c$ is initialized to a random value between 0
and 1. The subsequent temporal evolution of the density (which is
related to concentration via the EOS) is shown in Fig. \ref{fig:KHLM_3D},
showing the development of the instability with no visible numerical
artifacts. We also observe uniformly robust convergence of the GMRES
Stokes solver throughout the simulation.

\section{\label{sec:GiantFluct}Giant Concentration Fluctuations}

Advection of concentration by thermal velocity fluctuations in the
presence of large concentration gradients leads to the appearance
of \emph{giant fluctuations} of concentration, as has been studied
theoretically and experimentally for more than a decade \cite{GiantFluctuations_Nature,GiantFluctuations_Theory,GiantFluctuations_Cannell,FractalDiffusion_Microgravity,GiantFluctuations_Summary}.
In this Section, we use our algorithms to simulate experiments measuring
the temporal evolution of giant concentration fluctuations during
free diffusive mixing in a binary liquid mixture. In Ref. \cite{GiantFluctuations_Cannell},
Croccolo \emph{et al.} report experimental measurements of the temporal
evolution of the time-correlation functions of concentration fluctuations
during the diffusive mixing of water and glycerol. In the experiments,
a solution of glycerol in water with mass fraction of $c=0.39$ is
carefully injected in the bottom half of the experimental domain,
under the $c=0$ pure water in the top half. The two fluids slowly
mix over the course of many hours while a series of measurements of
the concentration fluctuations are performed.

In the experiments, quantitative shadowgraphy is used to observe and
measure the strength of the fluctuations in the concentration via
the change in the index of refraction. The observed light intensity,
once corrected for the optical transfer function of the equipment,
is proportional to the intensity of the fluctuations in the concentration
averaged along the vertical (gradient) direction, 
\[
c_{\perp}(x,z;\, t)=H^{-1}\int_{y=0}^{H}c(x,y,z;\, t)dy,
\]
where $H$ is the thickness of the sample in the vertical direction.
The quantity of interest is the correlation function of the Fourier
coefficients $\widehat{\d c}_{\perp}\left(k_{x},k_{z};\, t\right)$
of $c_{\perp}(x,z;\, t)$,
\[
C\left(\tau;\, t,k\right)=\av{\left(\widehat{\d c}_{\perp}\left(k_{x},k_{z};\, t+\tau\right)\right)\left(\widehat{\d c}_{\perp}\left(k_{x},k_{z};\, t\right)\right)^{\star}},
\]
where $k=\sqrt{k_{x}^{2}+k_{z}^{2}}$ is the wavenumber (in our two
dimensional simulations $k_{z}=0$), $\tau$ is a delay time, and
$t$ is the elapsed time since the beginning of the experiment. In
principle, the averaging above is an ensemble average but in the experimental
analysis, and also in our processing of the simulation results, a
time averaging over a specified time window $T$ is performed in lieu
of ensemble averaging. This approximation is justified because the
system is ergodic and the evolution of the deterministic (background)
state occurs via slow diffusive mixing of the water and glycerol solutions,
and thus happens on a much longer time scale (hours) than the time
delays of interest (a few seconds).

The Fourier transform (in time) of $C\left(\tau\right)$ is called
the dynamic structure factor. The equal-time correlation function
\[
S\left(k;\, t\right)=C\left(\tau=0;\, t,k\right)
\]
is the static structure factor, and is more difficult to measure in
experiments \cite{GiantFluctuations_Cannell}. For this reason, the
experimental results are presented in the form of normalized time-correlation
functions,
\[
\widetilde{C}\left(\tau;\, t,k\right)=\frac{C\left(\tau;\, t,k\right)}{S\left(k;\, t\right)}.
\]
The wavenumbers observed in the experiment and simulation are $k=\kappa\cdot2\pi/L$,
where $\kappa$ is an integer and $L$ is the horizontal extent of
the observation window or the simulation box size. When evaluating
the theory, we account for errors in the discrete approximation to
the continuum Laplacian by using the effective wavenumber
\begin{equation}
k_{\perp}=k_{x}\frac{\sin\left(k_{x}\Delta x/2\right)}{\left(k_{x}\Delta x/2\right)}\label{eq:modified_kx}
\end{equation}
instead of the actual discrete wavenumber $k_{x}$ \cite{LLNS_Staggered}.

The confinement in the vertical direction is expected to play a small
role because of the large thickness (2cm) of the sample, and a simple
quasi-periodic (bulk) approximation can be used. Approximate theoretical
analysis \cite{FluctHydroNonEq_Book} suggests that at steady state
the dominant nonequilibrium contribution to the static structure factor,
\begin{equation}
S\left(k;t\right)=\frac{k_{B}T}{\left(\eta\chi k^{4}-\rho\beta gh\right)}\, h^{2},\label{eq:S_c_c}
\end{equation}
exhibits a $k^{-4}$ power-law decay at large wavenumbers, and a plateau
to $k_{B}T\, h/\left(\rho\beta g\right)$ for wavenumbers smaller
than a rollover $k_{c}^{4}=\rho\beta gh/\left(\eta\chi\right)$ due
to the influence of gravity (buoyancy). Here $h\left(t\right)=d\bar{c}(y;t)/dy$
is the deterministic (background) concentration gradient, which decays
slowly with time due to the continued mixing of the water and glycerol
solutions.

An overdamped approximation suggests that the time correlations decay
exponentially, $\widetilde{C}\left(\tau;\, t,k\right)=\exp\left(-\tau/\tau_{k}\right)$,
with a relaxation time or decay time
\begin{equation}
\tau_{k}^{-1}=\chi k^{2}\left[1+\frac{\rho\beta gh}{\eta\chi k^{4}}\right],\label{eq:tau_overdamped}
\end{equation}
that has a minimum at $k=k_{c}$ with value $\tau_{\text{min}}^{-1}=2\chi k_{c}^{2}\sim\sqrt{hg}$.
For wavenumbers $k<k_{c}$ the relaxation time becomes \emph{smaller}
and can in fact become very small at the smallest wavenumbers, requiring
small time step sizes in the simulations to resolve the dynamics and
ensure stability of the temporal integrators. In the presence of gravity,
at small wavenumbers the separation of time scales used to justify
the overdamped limit fails and the fluid inertia has to be taken into
account \cite{MultiscaleIntegrators}. This changes the prediction
for the time correlation function to be a sum of two exponentials
with relaxation times, 
\begin{equation}
\tau_{1/2}^{-1}=\frac{1}{2}\left(\nu+\chi\right)k^{2}\pm\frac{1}{2}\,\sqrt{k^{4}\left(\nu-\chi\right)^{2}-4\beta gh},\label{eq:complex_tau}
\end{equation}
where $\nu=\eta/\rho$. This expression becomes complex-valued for
\[
k\lesssim k_{p}=\left(\frac{4\beta gh}{\nu^{2}}\right)^{\frac{1}{4}}=\left(4\frac{\chi}{\nu}\right)^{\frac{1}{4}}k_{c},
\]
indicating the appearance of \emph{propagative} rather than diffusive
modes for small wavenumbers, closely related to the more familiar
gravity waves. While experimental measurements over wavenumbers $k\lesssim k_{p}$
are not reported by Croccolo \emph{et al} \cite{GiantFluctuations_Cannell},
their experimental data does contain several wavenumbers in that range.
We report here simulation results for propagative concentration modes
at small wavenumbers. To our knowledge, experimental observation of
propagative modes has only been reported for temperature fluctuations
\cite{TemperatureGradient_Cannell}.

Because it is essentially impossible to analytically solve the linearized
fluctuating equations in the presence of spatially-inhomogeneous density
and transport cofficients and nontrivial boundary conditions, the
existing theoretical analysis of the diffusive mixing process \cite{GiantFluctuations_Theory}
makes a quasi-periodic constant-coefficient and constant-gradient
incompressible approximation \cite{FluctHydroNonEq_Book}. This approximation,
while sufficient for qualitative studies, is not appropriate for quantitative
studies because the viscosity $\eta$ and mass diffusion coefficient
$\chi$ vary by about a factor of three from the bottom to the top
of the sample. In our simulations we account for the full dependence
of density, viscosity and diffusion coefficient on concentration.

\subsection{Simulation Parameters}

For LFHD there is no difference between the two and three dimensional
problems due to the symmetries of the problem \cite{FluctHydroNonEq_Book}.
Because very long simulations with a small time step size are required
for this study, we perform two-dimensional simulations. Furthermore,
in these simulations we do not include a stochastic flux in the concentration
equation, i.e., we set $\M{\Psi}=0$, so that all fluctuations in
concentration arise from the coupling to the fluctuating velocity.
With this approximation we do not need to model the chemical potential
of the mixture and obtain $\mu_{c}$. This approximation is justified
by the fact that it is known experimentally that the nonequilibrium
fluctuations are much larger than the equilibrium ones for the conditions
we consider \cite{GiantFluctuations_Cannell}; in fact, the fluctuations
due to nonzero $\M{\Psi}$ are smaller than solver or even roundoff
tolerances in the simulations reported here.

We base our parameters on the experimental studies of diffusive mixing
in a water-glycerol mixture, as reported by Croccolo \emph{et al.}
\cite{GiantFluctuations_Cannell}. In the actual experiments the fluid
sample is confined in a cylinder $2\,\mbox{cm}$ in diameter and $2\,\mbox{cm}$
thick in the vertical direction. In our simulations, the two-dimensional
physical domain is $1.132\,\textrm{cm}\times1.132\,\textrm{cm}$ discretized
on a uniform $256\times256$ two dimensional grid, with a thickness
of $1\,\textrm{cm}$ along the $z$ direction. This large thickness
makes the magnitude of the fluctuations very small since the cell
volume $\D V$ contains a very large number of molecules, and puts
us in the linearized regime \cite{MultiscaleIntegrators}. The width
of the domain $L=1.132\,\textrm{cm}$ is chosen to match the observation
window in the experiments, and thus also match the discrete set of
wavenumbers between the simulations and experiments. Earth gravity
$g=-9.81\,\mbox{m}^{2}/\mbox{s}$ is applied in the negative $y$
(vertical) direction; for comparison we also perform a set of simulations
without gravity. Periodic boundary conditions are applied in the $x$-direction
and impenetrable no-slip walls are placed at the $y$ boundaries.
The initial condition is $c=0.39$ in the bottom half of the domain
and $c=0$ in the top half, with velocity initialized to zero. The
temperature is kept constant at $300\,\mbox{K}$ throughout the domain.
Centered advection is used to ensure fluctuation-dissipation balance
over the whole range of wavenumbers represented on the grid. 

A very good fit to the experimental equation of state (dependence
of density on concentration at standard temperature and pressure)
over the whole range of concentrations of interest is provided by
the EOS (\ref{eq:EOS_quasi_incomp}) with the density of water set
to $\bar{\rho}_{2}=1\,\textrm{g}/\textrm{cm}^{3}$ and the density
of glycerol set to $\bar{\rho}_{1}=1.29\,\textrm{g}/\textrm{cm}^{3}$.
Experimentally, the dependence of viscosity on glycerol mass fraction
has been fit to an exponential function \cite{GiantFluctuations_Cannell},
which we approximate with a rational function over the range of concentrations
of interest \cite{WaterGlycerolViscosity},
\begin{equation}
\eta(c)\approx\frac{1.009+1.1262\, c}{1-1.5326\, c}\cdot10^{-2}\,\frac{\textrm{g}}{\textrm{cm}\,\mbox{s}}.\label{eq:nu_c}
\end{equation}
The diffusion coefficient dependence on the concentration has been
studied experimentally, and we employ the fit proposed in Ref. \cite{WaterGlycerolDiffusion},
\begin{equation}
\chi(c)=\frac{1.024-1.002\, c}{1+0.663\, c}\cdot10^{-5}\,\frac{\textrm{cm}^{2}}{\mbox{s}},\label{eq:chi_c}
\end{equation}
which is in reasonable but not perfect agreeement with a Stokes-Einstein
relation $\eta(c)\chi(c)=\mbox{const.}$ Note that the Schmidt number
$S_{c}=\nu/\chi\sim10^{3}$. In Ref. \cite{GiantFluctuations_Cannell},
based on the experimental measurements and the approximate theoretical
model, it is suggested that $\chi\approx10^{-5}\,\textrm{cm}^{2}/\textrm{s}$
is constant over the range of concentrations present. For comparison,
we also perform simulations in which we keep the diffusion coefficient
independent of concentration, while still taking into account the
concentration dependence of viscosity. It is worth noting that there
is a notable disagreement between experimental measurements of $\chi(c)$
using different experimental techniques \cite{WaterGlycerolDiffusion}
and the true dependence is not known with the same accuracy as that
of $\eta(c)$.

When gravity is present, we use the inertial Algorithm \ref{alg:LMInertial},
with a rather small time step size $\Delta t=0.01375$s due to the
fact that the smallest relaxation time measured is on the order of
$0.1$s. For this time step size, the viscous CFL number is $\nu\D t/\D x^{2}\sim10-30$,
indicating that the viscous dynamics is resolved except at the wavenumbers
comparable to the grid spacing. In the absence of gravity we use the
overdamped Algorithm \ref{alg:LMOverdamped}, which allows us to use
a much larger time step size (on the diffusive time scale), $\D t=0.22$s,
giving a diffusive CFL number on the order of $\chi\D t/\D x^{2}\lesssim0.1$.
Using larger time step sizes than this would require an implicit treatment
of mass diffusion.

\subsection{Results}

Our simulations closely mimic the experiments of Croccolo \emph{et
al} \cite{GiantFluctuations_Cannell}. We perform a long (stochastic)
run of the diffusive mixing up to physical time $t=21,021$s, saving
a snapshot and statistics every 21,840 time steps, which corresponds
to 300 seconds of physical time. We then perform 8 short runs with
different random seeds starting from the saved snapshots, and compute
time correlation functions averaged over a short time interval. Note
that in the experiments a similar procedure is used in which data
is collected over short time intervals during a single long mixing
process. Croccolo \emph{et al.} report measurements at $t=600$s,
$3060$s, $8160$s, and $14880$s. Table \ref{tab:WaterGly} lists
the time intervals over which we collect statistics in the simulations,
which match those in the experiments as well as possible. The time
interval between successive snapshots used in the computation of the
time correlation function is four time steps or 0.055s, which is four
times smaller than the interval used in the experimental analysis.
In the experiments averaging is performed over a range of wavenumbers
in the $\left(k_{x},k_{z}\right)$ plane with similar magnitude. Since
we perform two dimensional simulations we average over the 8 independent
simulations; in the end the statistical errors are lower in the simulation
results since experiment are subject to a large experimental noise
not present in the simulations.

\begin{table}
\begin{centering}
\begin{tabular}{l|c|r}
\hline 
Starting Time  &
Total Time Steps  &
End Time \tabularnewline
\hline 
600.6 s  &
3328  &
$600+3328\Delta t=646.36$ s \tabularnewline
\hline 
3003 s  &
10784  &
3151.28 s \tabularnewline
\hline 
8108.1 s  &
4992  &
8176.74 s \tabularnewline
\hline 
15015 s  &
4768  &
15080.6 s\tabularnewline
\hline 
\end{tabular}
\par\end{centering}

\caption{\label{tab:WaterGly}Time intervals over which we average the dynamic
correlation functions used to compute the relaxation times shown in
Fig. \ref{fig:WaterGly_tau}.}
\end{table}

\subsubsection{Dynamic Structure Factors}

\begin{figure}
\begin{centering}
\includegraphics[width=0.49\textwidth]{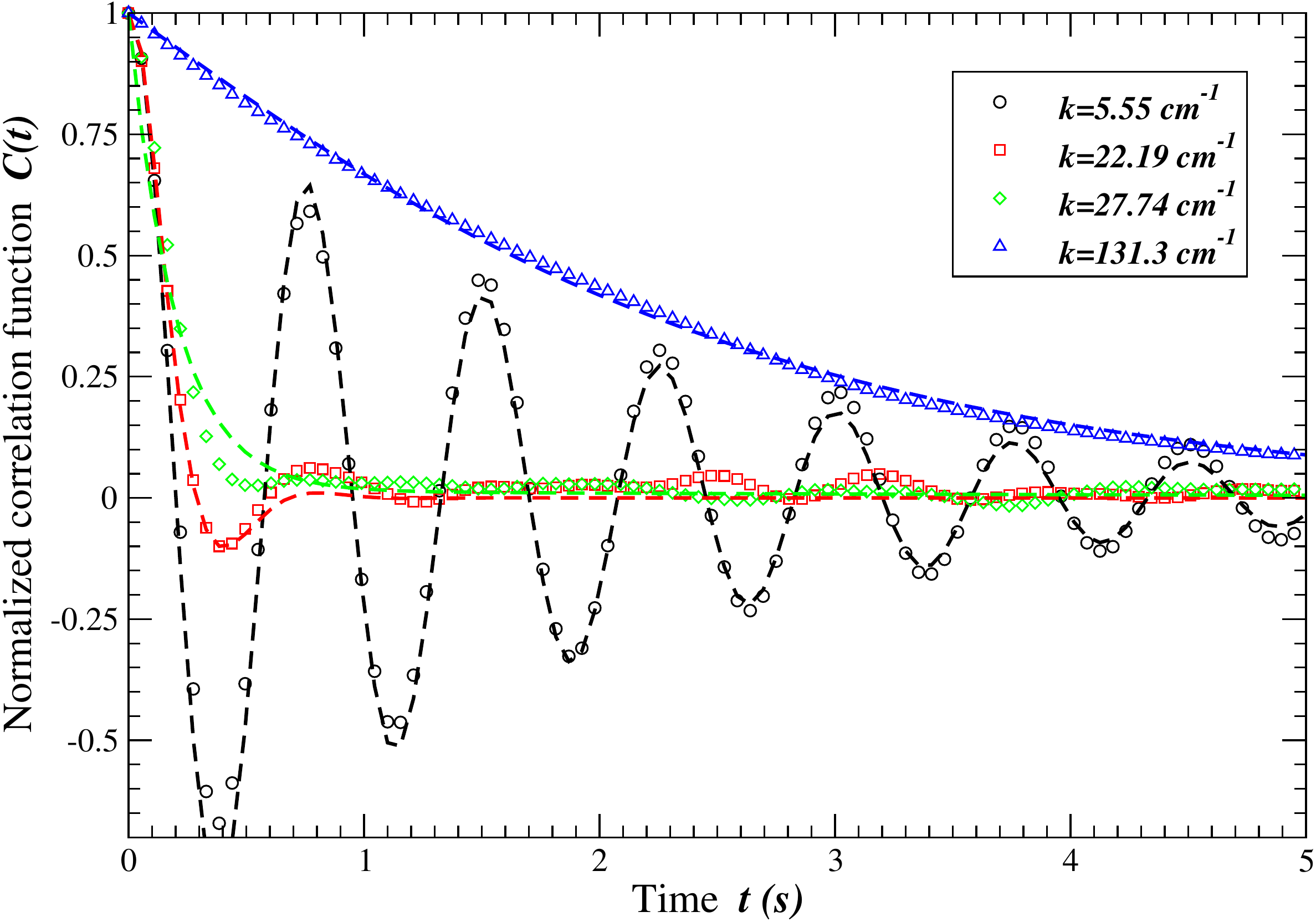}\includegraphics[width=0.49\textwidth]{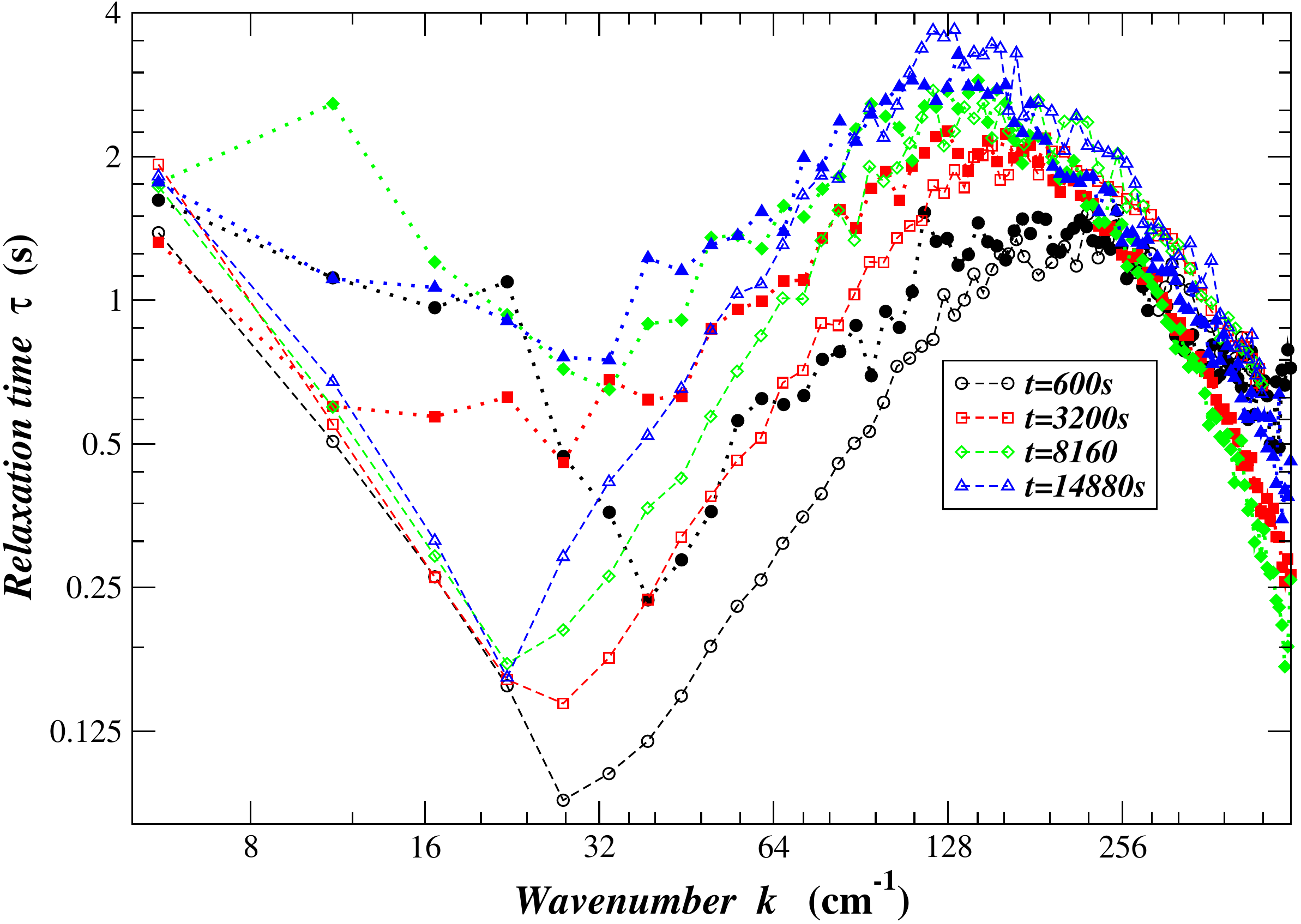}
\par\end{centering}

\caption{\label{fig:WaterGly_tau}Dynamics of concentration fluctuations during
free-diffusive mixing of water and glycerol. (Left) Numerical results
for the time correlation functions for several selected wavenumbers
about $8160$s from the beginning of the experiment. Symbols indicate
results from the simulations and lines of the same color indicate
the fit to (\ref{trig}) for the first ($k\approx5.6\,\mbox{cm}^{-1}$)
and fourth ($k\approx22.2\,\mbox{cm}^{-1}$) wavenumbers, or to (\ref{double})
for the remaining wavenumbers. Note that the statistical errors due
to the finite averaging increase with time and the tails of the correlation
functions are not reliably estimated. (Right) Relaxation or decay
times as a function of wavenumber at several points in time. Empty
symbols show results from computer simulations, and filled from experimental
measurements \cite{GiantFluctuations_Cannell}.}
\end{figure}

In order to extract a relaxation time, we fit the numerical results
for the normalized time-correlation function to an analytical formula.
For the first four wavenumbers $k=\left(1,2,3,4\right)\cdot2\pi/L$,
clear oscillations (propagative modes) were observed, as illustrated
in the left panel of Fig. \ref{fig:WaterGly_tau}. For these wavenumbers
we used the fit
\begin{equation}
\widetilde{C}\left(t\right)=\exp\left(-\frac{t}{\tau}\right)\left(A\sin\left(2\pi\frac{t}{T}\right)+\cos\left(2\pi\frac{t}{T}\right)\right),\label{trig}
\end{equation}
where the relaxation time $\tau$, the coefficient $A$ and the period
of oscillation $T$ are the fitting parameters. For the remaining
wavenumbers, we used a double-exponential decay for the fitting,
\begin{equation}
\widetilde{C}\left(t\right)=\alpha\exp\left(-\frac{t}{\tau_{1}}\right)+(1-\alpha)\exp\left(-\frac{t}{\tau_{2}}\right),\label{double}
\end{equation}
where $\alpha$, $\tau_{1}$ and $\tau_{2}$ are the fitting parameters.
This leads to good fits for $k>k_{p}\sim32\,\mbox{cm}^{-1}$; for
a few transitional wavenumbers such as $k\sim28\,\mbox{cm}^{-1}$
the fit is not as good, as illustrated in the left panel of Fig. \ref{fig:WaterGly_tau}.
From the fit (\ref{double}) we obtain the relaxation time $\tau$
as the point at which the amplitude decays by $\widetilde{C}\left(\tau\right)=1/e$.

A similar procedure was also used to obtain the relaxation time from
the experimental data of Croccolo et al. \cite{GiantFluctuations_Cannell}
for all wavenumbers %
\footnote{The experimental data for the time correlation functions were graciously
given to us by Fabrizio Croccolo.%
}. The experimental data shows monotonically decaying correlation functions
$\widetilde{C}\left(\tau;\, t,k\right)$ for \emph{all} measured wavenumbers,
not consistent with the oscillatory correlation function observed
for the four smallest wavenumbers in our simulations \cite{MultiscaleIntegrators},
see the left panel of Fig. \ref{fig:WaterGly_tau}. We believe that
this mismatch is due to the way measurements for different wavenumbers
of similar modulus are averaged in the experimental calculations.
In our two-dimensional simulations, we do not perform any averaging
over wavenumbers. We believe that the experimentally measured time
correlation functions capture the \emph{real} part of the decay times
\emph{only} and thus have the form of a sum of exponentials. Due to
the lower time resolution and the fact that the static structure factor
is not known, for the experimental data we used a single exponential
fit and added an offset to account for the background noise,
\[
C\left(t\right)=A\exp\left(-\frac{t}{\tau}\right)+B.
\]

In the right panel of Fig. \ref{fig:WaterGly_tau} we compare simulation
and experimental results for the real part of the decay or relaxation
time $\tau_{k}$, at several points in time measured from the beginning
of the experiment. Good agreement is observed between the two with
the same qualitative trends: A diffusive relaxation time $\tau_{k}^{-1}\approx\chi k^{2}$
for large wavenumbers, with a maximum at $k\approx k_{c}$, and then
another minimum at $k\approx k_{p}$. Note that decay times are not
reported by Croccolo \emph{et al} \cite{GiantFluctuations_Cannell}
for wavenumbers $k\lesssim k_{p}$ since that work focuses on the
effect of gravity for $k\lesssim k_{c}$. In our analysis of the experimental
data we included all measured wavenumbers, including those for which
propagative modes are observed. Here, the diffusion coefficient varies
with concentration according to (\ref{eq:chi_c}); very similar results
for the relaxation times were obtained by keeping $\chi\approx10^{-5}\,\textrm{cm}^{2}/\textrm{s}$
constant, as suggested by Croccolo \emph{et al} \cite{GiantFluctuations_Cannell}.
This indicates that the dynamic correlations are not very sensitive
to the concentration dependence $\chi(c)$. In future work we will
perform a more careful comparison to experiments.

\subsubsection{Static Structure Factors}

\begin{figure}
\begin{centering}
\includegraphics[width=0.49\textwidth]{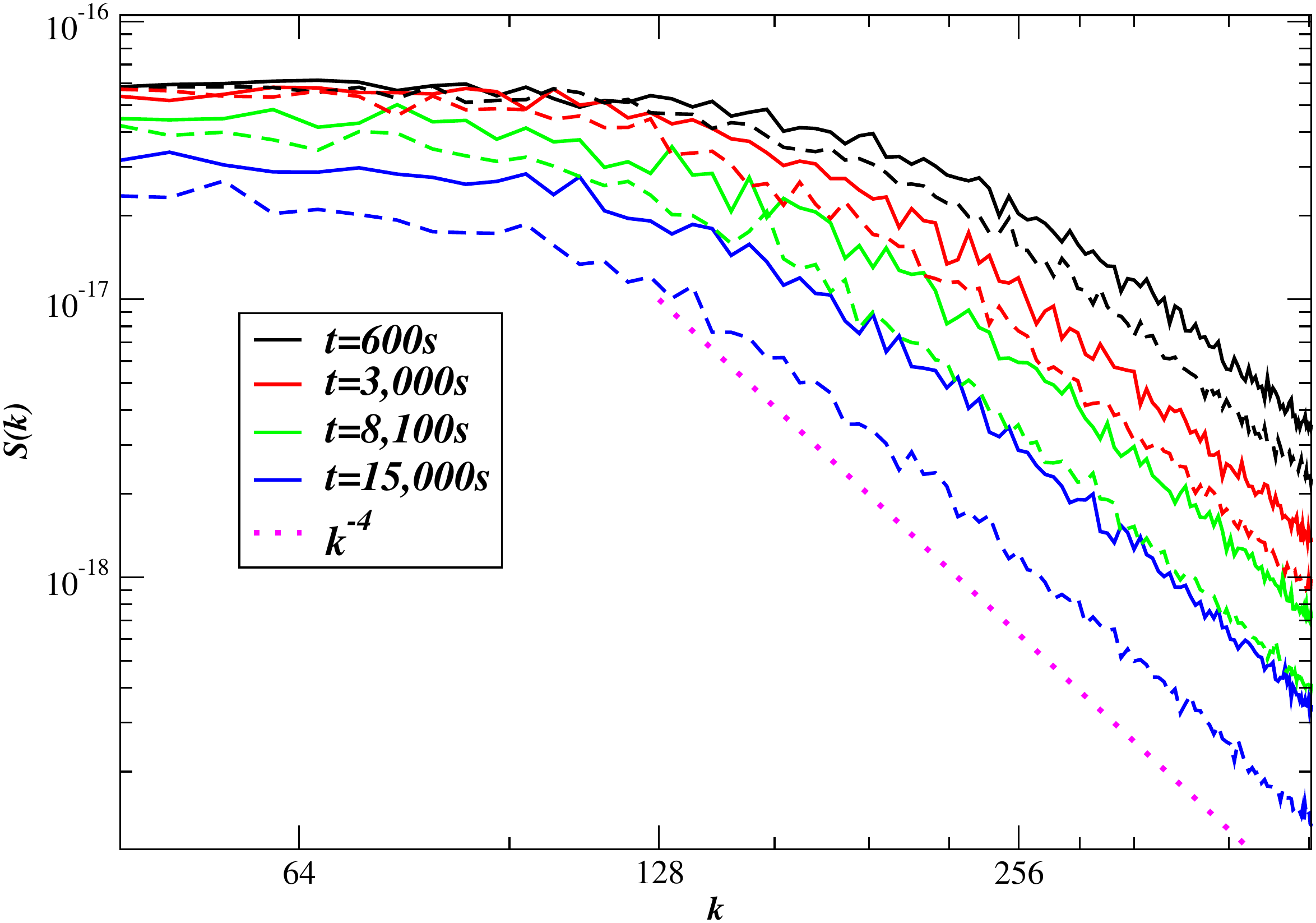}\includegraphics[width=0.49\textwidth]{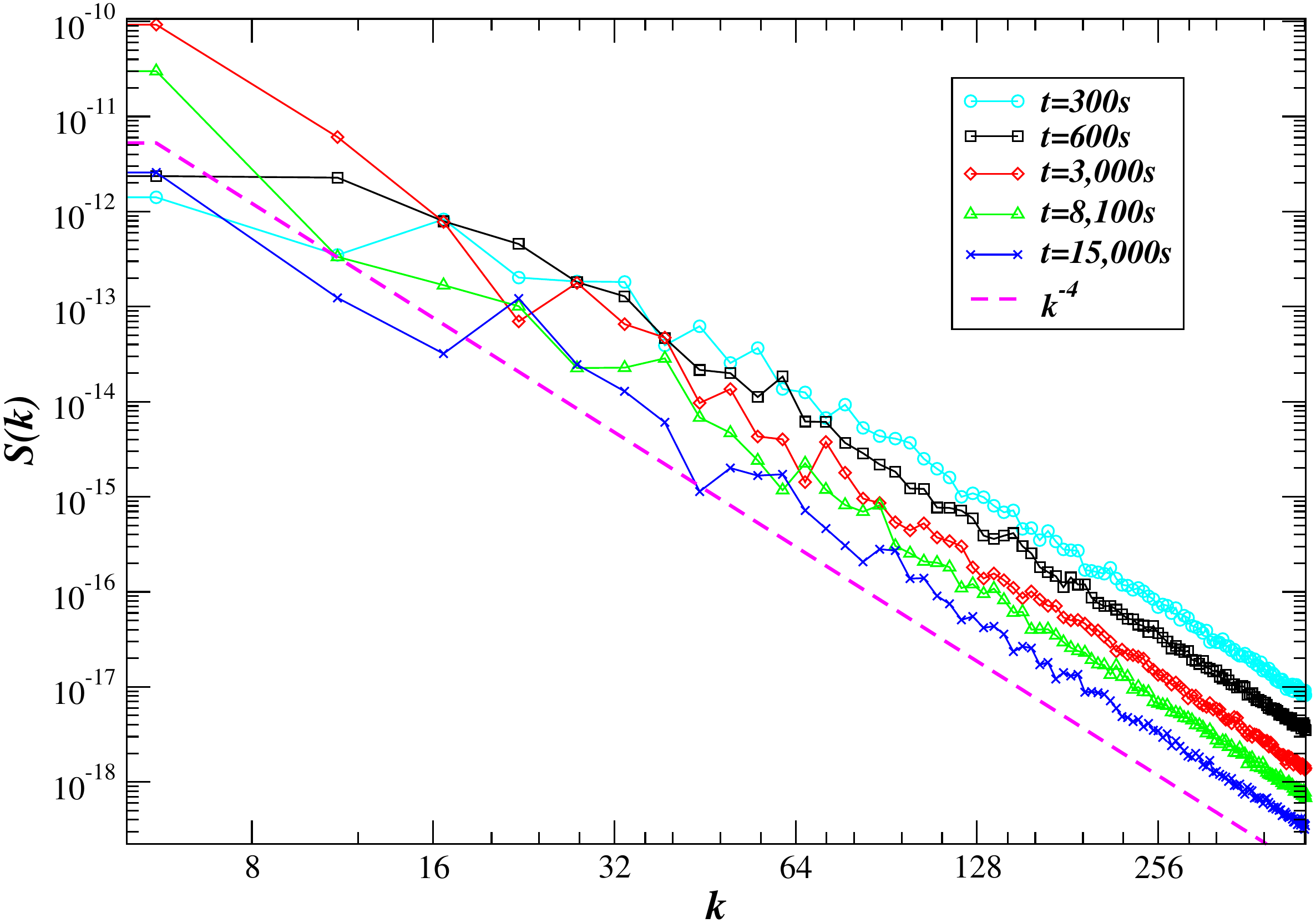}
\par\end{centering}

\caption{\label{fig:WaterGly_S_k}Evolution of the static structure factor
during free-diffusive mixing of water and glycerol. (Left) With Earth
gravity. Solid lines show results for simulations in which $\chi(c)$
depends on concentration according to (\ref{eq:chi_c}), while dashed
lines of the same color correspond to simulations in which $\chi\approx10^{-5}\,\textrm{cm}^{2}/\textrm{s}$
is constant. Fluctuations at large wavenumbers follow a $k^{-4}$
power law but are damped by gravity at small wavenumbers. (Right)
Without gravity; observe the large difference in the vertical axes
showing ``giant'' fluctuations in the microgravity case. Note that
these are results from a single simulation, mimicking a single experimeent,
and therefore there are large statistical uncertainties at small wavenumbers
(large decorrelation times).}
\end{figure}

Extracting the static structure factor from experimental measurements
is complicated by several factors, including the presence of optical
prefactors such as the transfer function of the instrument, and the
appearance of additional contributions to the scattered light intensity
such as shot noise, contributions due to giant temperature fluctuations
\cite{TemperatureGradient_Cannell}, and capillary waves \cite{GiantFluctuations_Interfaces,GiantFluctuations_Capillary}.
We therefore study the evolution of the static structure factor using
simulations only. In the left panel of Fig. \ref{fig:WaterGly_S_k}
we show numerical results for the static structure factor $S\left(k;\, t\right)$
of the discrete concentration field averaged along the $y$-axes,
at a series of times $t$ chosen to match those of the experimental
measurements. Instead of ensemble averaging, here we performed a temporal
average of the spectrum of the vertically-averaged concentration over
a period of $300$s, \emph{ending} at the time indicated in the legend
of the figure. The characteristic $k^{-4}$ power law decay at large
wavenumbers and the plateau at small wavenumbers predicted by (\ref{eq:S_c_c})
are clearly observed in Fig. \ref{fig:WaterGly_S_k}, consistent with
a value of $h$ decreasing with time. A quantitative difference is
seen between the results for variable and constant diffusion coefficients,
consistent with a different value of the imposed concentration gradient
$h$ due to the somewhat different evolution of $\bar{c}(y,t)$.

In the right panel of Fig. \ref{fig:WaterGly_S_k} we show numerical
results for the static structure factor that would be obtained had
the experiment been performed in microgravity ($g=0$). In this case,
we use the overdamped algorithm \ref{alg:LMOverdamped} since there
is a persistent large separation of time scales between the slow concentration
and fast velocity. We see clear development of a $k^{-4}$ power law
as predicted by (\ref{eq:S_c_c}) for $g=0$. Note that here the concentration
gradient is established instantaneously, in fact, it is the largest
in the initial configuration and then decays on the diffusive time
scale; this is different from simulations of the development of giant
fluctuations in microgravity during the GRADFLEX experiment reported
in \cite{MultiscaleIntegrators}, in which the concentration gradient
is slowly established on a diffusive scale. The results in Fig. \ref{fig:WaterGly_S_k}
show that it takes some time for the giant fluctuations at smallest
wavenumbers to develop; the diffusive relaxation time corresponding
to the smallest wavenumber studied, $k_{\min}\approx5\,\mbox{cm}^{-1}$
is $\tau_{\max}=\left(\chi k_{\min}^{2}\right)^{-1}\sim4,000$s. After
a time $\sim\left(\chi k^{2}\right)^{-1}$, the amplitude of the fluctuations
$S(k)\sim k^{-4}h^{2}(t)$ decays slowly due to the diffusive mixing,
and eventually the system will fully mix and reach thermodynamic equilibrium.

\section{\label{sec:Conclusions}Conclusions}

We have developed a low Mach number algorithm for diffusively-mixing
mixtures of two liquids with potentially different density and transport
coefficients. In the low Mach number setting, the incompressible constraint
is replaced by a quasi-compressible constraint that dictates that
stochastic and diffusive mass fluxes must create local expansion and
contraction of the fluid to maintain a constant thermodynamic (base)
pressure.

We employed a uniform-grid staggered-grid spatial discretization \cite{LLNS_Staggered}.
Following prior work in the incompressible simple-liquid case \cite{NonProjection_Griffith},
we treated viscosity implicitly without splitting the pressure update,
relying on a recently-developed variable-coefficient Stokes solver
\cite{StokesKrylov} for efficiency. This approach works well for
any Reynolds number, including the viscous-dominated overdamped (zero
Reynolds number) limit, even in the presence of nontrivial boundary
conditions. Furthermore, by using a high-resolution BDS scheme \cite{BDS}
to advect the concentration we robustly handled the case of no mass
diffusion (no dissipation in the concentration equation). In our spatial
discretization we strictly preserved mass and momentum conservation,
as well as the equation of state (EOS) constraint, by using a finite-volume
(flux-based) discretization of advective fluxes in which fluxes are
computing using extrapolated values of concentration and density that
obey the EOS. Our temporal discretization used a predictor-corrector
integrator that treats all terms except momentum diffusion (viscosity)
explicitly \cite{MultiscaleIntegrators}.

We empirically verified second-order spatio-temporal accuracy in the
deterministic method. In the stochastic context, establishing the
weak order of accuracy is nontrivial in the general low Mach number
setting. For centered advection our temporal integration schemes can
be shown to be second-order accurate for the special case of a Boussinesq
constant-density (incompressible) approximation, or in the overdamped
(inertia-free) limit. Existing stochastic analysis does not apply
to the case of BDS advection because Godunov schemes do not fit a
method-of-lines approach, but rather, employ a space-time construction
of the fluxes. The presence of nontrivial density differences between
the pure fluid components and nonzero mass diffusion coefficient,
complicates the analysis even for centered advection, due to the presence
of a nontrivial EOS constraint on the fluid dynamics. It is a challenge
for future work to develop improved numerical analysis of our schemes
in both the deterministic and the stochastic setting.

In future work, we will demonstrate how to extend the algorithms proposed
here to multispecies mixtures of liquids using a generalization of
the low Mach number constraint. The nontrivial multispecies formulation
of the diffusive and stochastic mass fluxes has already been developed
by some of us in the compressible setting \cite{MultispeciesCompressible}.

It is also possible to include thermal effects in our formulation,
by treating the temperature in a manner similar to the way we treated
concentration here. Two key difficulties are constructing a spatial
discretization that ensures preservation of an appropriately generalized
EOS, as well as developing temporal integrators that can handle the
moderate separation of time scales between the (typically) slower
heat diffusion and (typically) faster momentum diffusion. In particular,
it seems desirable to also treat temperature implicitly. Such implicit
treatment of mass or heat diffusion is nontrivial because it would
require solving coupled (via the EOS constraint) velocity-temperature
or velocity-concentration linear systems, and requires further investigation.

In the staggered-grid based discretization developed here, we can
only employ existing higher-order Godunov advection schemes for the
cell-centered scalar fields such as concentration and density. It
is a challenge for future work to develop comparable methods to handle
advection of the staggered momentum field. This would enable simulations
of large Reynolds number flows. It should be noted, however, that
our unsplit approach is most advantageous at small Reynolds numbers. 

A challenge for future work on low Mach number fluctuating hydrodynamics
is to account for the effects of surface tension in mixtures of immiscible
or partially miscible liquids. This can be most straightforwardly
accomplished by using a diffuse-interface model, as some of us recently
did in the compressible setting for a single-fluid multi-phase system
\cite{CHN_Compressible}. One of the key challenges is handling the
fourth-order derivative term in the concentration equation in a way
that ensures stability of the temporal integrator, as well as developing
a consistent discretization of the Korteweg stresses on a staggered
grid \cite{StagerredFluct_Inhomogeneous}. 

The semi-implicit temporal integrators we described here can deal
well with a broad range of Reynolds or Schmidt numbers in the deterministic
(smooth) setting. In the context of fluctuating hydrodynamics, however,
all modes are thermally excited and treatment of viscosity based on
a Crank-Nicolson method (implicit midpoint rule) are bound to fail
for sufficiently large Schmidt numbers (or sufficiently low Reynolds
numbers). In this work we solved this problem for the case of infinite
Schmidt, zero Reynolds number flows by taking an overdamped limit
of the original inertial equations before temporal discretization.
It is a notable challenge for the future to develop uniformly accurate
temporal integrators that work over a broad range of Reynolds or Schmidt
numbers, including the asymptotic overdamped limit, in the presence
of thermal fluctuations.
\begin{acknowledgments}
We would like to thank Fabrizio Croccolo and Alberto Vailati for sharing
their experimental data on water-glycerol mixing, as well as numerous
informative discussions. This material is based upon work supported
by the U.S. Department of Energy Office of Science, Office of Advanced
Scientific Computing Research, Applied Mathematics program under Award
Number DE-SC0008271 and under contract No. DE-AC02-05CH11231. Additional
support for A. Donev was provided by the National Science Foundation
under grant DMS-1115341.
\end{acknowledgments}


\end{document}